\documentclass[preprint,12pt]{elsarticle}

\usepackage{amssymb}

\usepackage{graphicx}
\usepackage{subfigure}
\usepackage{amsfonts}
\usepackage{amsmath}
\usepackage{amssymb}
\usepackage{psfig}

\newcommand \nn \nonumber
\newtheorem {thm}{Theorem}

\newtheorem {prop}{Proposition}

\journal{Information Sciences}

\begin{document}

\begin{frontmatter}

\title{Performance of a Class of Multi-Robot Deploy and Search Strategies based on Centroidal Voronoi Configurations}

\author[1]{K.R. Guruprasad}
\ead{krgprao@gmail.com}
\author[2]{Debsish Ghose}
\ead{dghose@aero.iisc.ernet.in}

\address[1]{Department of Mechanical Engineering\\
National Institute of Technology Karnataka\\
Surathkal, INDIA
}
\address[2]{Guidance, Control, and Decision Systems Laboratory\\
Department of Aerospace Engineering\\
Indian Institute of Science, Bangalore, INDIA}

\begin{abstract}
This paper considers a class of deploy and search strategies for
multi-robot systems and evaluates their performance.  The application framework
used is a system of autonomous mobile robots equipped with required
sensors and communication equipment deployed in a search space to
gather information. The lack of information about the search space
is modeled as an uncertainty density distribution over the search
space. A {\em combined deploy and search} (CDS) strategy has been formulated as a modification to {\em sequential deploy and search} (SDS) strategy presented in our previous work. The optimal deployment strategy using Voronoi partition forms the basis for these two search strategies. The strategies are analyzed in presence of constraints on robot speed and limit on sensor range for convergence of trajectories with corresponding control laws responsible for the motion of robots. SDS and CDS strategies are compared with standard greedy and random search strategies on the basis of time taken to achieve reduction in the uncertainty density below a desired
level. The simulation experiments reveal several important issues related
to the dependence of  the relative performances
of the strategies  on  parameters such as number of robots,
speed of robots, and their sensor range limits.
\end{abstract}


\end{frontmatter}

\section{Introduction}
\label{intro}
In nature, we observe groups of animals performing a large number of
complex tasks by cooperating among themselves in a distributed
manner. Each member of the group performs a relatively  simple task
using information available in its neighborhood leading to
emergence of a collective behavior. Flocks of birds, schools of
fish, foraging of ants are a few examples of this phenomena.
Nature seems to solve complicated problems in
simple and elegant ways. Instead of increasing the
complexity/intelligence of an individual to solve a complex problem,
nature relies on a group of individuals with nominal abilities. The outcome of such cooperation is more
robust to failure of a few of the group members. From a more
philosophical perspective, the focus is on a society rather than on an individual member.

The problem of searching for targets in unknown environments has been addressed in the literature in the past
\cite{koopman},\cite{stone},\cite{benkoski}, and \cite{lida}. Various search strategies available in the literature have been surveyed in \cite{benkoski}. These fundamental works
were mostly theoretical in nature and were applicable to a single agent searching for single or multiple, static or
moving, targets. It is likely that the same task can be accomplished more effectively by multiple searchers. But
when multiple agents are involved, coordination between them becomes an important issue. Although a centralized
controller can solve the problem, it has many shortcomings such as requirement of complete connectivity, large
communication effort, etc. Further, failure of the central controller leads to failure of the entire system. As discussed earlier, most biological systems such as ants, birds, fish etc., have distributed local decision making capabilities which, in turn, lead to a useful collective behavior such as swarms, flocks, schools, etc. Agents taking decisions based only on available local information and distributed
control law (usually referred to as `behavior' in biological systems), lead to coordination among the agents and
result in a meaningful collective behavior. Developments in areas such as wireless communication, autonomous
vehicular technology, computation, and sensors, facilitate the use of large number of agents (UAVs, mobile robots, or
autonomous vehicles), equipped with sensors, communication equipment, and computation ability, to cooperatively
achieve various tasks in a distributed manner.

Distributed multi-agent systems have been shown to achieve and maintain formations, move as flocks while
avoiding obstacles, etc., thus mimicking their biological counterparts. They can also be used in applications such as
search and rescue, surveillance, multiple source identification, and cooperative transportation. The major advantages
of distributed systems are immunity to failure of individual agents, their versatility in accomplishing multiple tasks,
simplicity of agents' hardware, and requirement of only minimal local information. At the same time it is important
to design distributed control laws that guarantee stability and convergence to the desired collective behavior under
limited information and evolving network configurations.

\subsection{Related Literature}

One class of problems discussed in the literature is that, of optimally locating agents or sensors in the domain
of interest and belongs to the class of locational optimization or facility location problems \cite{drezner,okabe}. A centroidal
Voronoi configuration is a standard solution for this class of problems \cite{du}, where each of the
agents is located at the centroid of the corresponding Voronoi cell. Cortes et al. \cite{bullo1,bullo2} use these concepts to solve a
spatially distributed optimal deployment problem for multi-agent systems. A density distribution, as a measure of
the probability of occurrence of an event, is used, along with a function of the Euclidean
distance providing a quantitative assessment of how poor the sensing performance is, to formulate the problem. Centroidal Voronoi configuration, with centroids of Voronoi cells computed
based on the density distribution within the cell, is shown to be the optimal deployment of sensors minimizing
the sensory error. The Voronoi partition becomes natural optimal partitioning due to monotonic variation of sensor
effectiveness function with the Euclidean distance. Schwager et al. \cite{schwager} interpret the density distribution of \cite{bullo1} in a
non-probabilistic framework and approximate it by sensor measurements. Sujit and Beard \cite{sujit} present an exploration system for multiple unmanned aerial vehicles (UAVs) navigating through a simulated unknown region that contains obstacles of unknown shape, size, and initial position. They perform Monte-Carlo simulation to analyze the effect on area coverage with changes in number of agents, sensor range, and communication range.

Cooperative search by multiple agents has been studied by various researchers. Enns et al. \cite{enns} use predefined
lanes prioritizing them with the target probability. The vehicles cooperate to ensure that the total path length covered
by them is minimized while exhaustively searching the area. A dynamic inversion based control law is used to
make vehicles track the assigned tracks or lanes while considering the maximum turn radius constraint. Spires and
Goldsmith \cite{spires} use space filling curves such as Hilbert curves to cover a given space and perform exhaustive search
by multiple robots. Vincent and Rubin \cite{vincent} address the problem of cooperative search strategies for UAVs searching for moving, possibly evading, targets in a hazardous environment. They use predefined
swarm patterns with an objective of maximizing the target detection probability in minimum expected time and using
minimum number of UAVs having limited communication range. Beard and McLain \cite{beard} use dynamic programming
methods to develop strategies for a team of cooperating UAVs to visit as many opportunities without collision while
avoiding hazards, in a search area which contains regions of opportunity and hazards. The UAVs have the additional
requirement that they should stay within communication range of each other. Flint et al. \cite{flint} provide a model and
algorithm for path planning of multiple UAVs searching in an uncertain and risky environment using dynamic
programming approach. The search area is divided into cells and in each cell the probability of existence of a
target defined. Pfister \cite{pf} uses fuzzy cognitive map to model the cooperative control process in an autonomous
vehicle. In \cite{poly},\cite{yang1}, and \cite{yang2} the authors use distributed reinforcement learning and planning for cooperative multi-agent
search. The agents learn the environment online and store the information in the form of a search map and utilize
this information to compute their trajectory. The agents are assumed to have limited maneuvering ability, sensor
range and fuel. In \cite{yang1} the authors show a finite lower bound on the search time. Rajnarayan and Ghose \cite{rajghose} use
concepts from team theory to formulate the multi-agent search problem as a nonlinear programming problem in a centralized perfect information case. The problem is then reformulated in a Linear-Quadratic-Gaussian setting that
admits a decentralized team theoretic solution. Dell et al. \cite{dell} develop an optimal branch-and-bound procedure with
heuristics such as combinatorial optimization, genetic algorithm and local start with random restarts, for solving
constrained-path problems with multiple searchers. Jin et al. \cite{jin}
address a search and destroy problem in a military setting with heterogeneous team of UAVs.
Altshuler et al. \cite{altshuler} examine the Cooperative Hunters problem, where a swarm of UAVs is used for searching one or more "evading targets", which are moving in a predefined area while trying to avoid a detection by the swarm. Sathyaraj et al. \cite{jain} provide a comparative study of path planning algorithm for multiple UAVs used in team for detecting targets and keeping them in its sensor range. Yatsenko et al. \cite{yatsenko} discuss problems dealing with cooperative control of multiple agents moving in a region searching for targets.

In  \cite{tase} we have proposed a search strategy namely {\em sequential deploy and search} (SDS) for multiple agents such as UAVs or mobile robots using Voronoi partition with some preliminary results. In this work, an uncertainty density was used to model the lack of information about the search area. Search agents are deployed optimally maximizing a single step search effectiveness and then perform search. Each agent performs search within its Voronoi cell where it is most effective. The optimal deployment was shown to be a variation of centroidal Voronoi configuration, where each agent is located at the centroid of its Voronoi cell with perceived uncertainty density. The material presented in this paper is based on the doctoral thesis \cite{mythesis}.

In this paper, we provide a more detailed account of the SDS strategy and propose a modification of the this strategy, named as {\em combined deploy and search} (CDS) where, the mobile robots perform search in discrete steps, as they move toward the centroids. We compare the performance of SDS and CDS strategies with standard {\em greedy} and {\em random} search strategies based on simulation experiments.

\subsection{Contribution of the paper} We address a multi-robot
search problem where $N$ agents, equipped with sensors, search the
space $Q$, a convex polytope in $\mathbb{R}^d$, the $d$-dimensional Euclidean space. The lack of
information about the search space is modeled as an uncertainty
density distribution $\phi : Q \mapsto [0,1]$. We denote $P(t) =
\{p_1,p_2,\ldots,p_N\}$, $p_i \in Q$, with $p_i \neq p_j$, whenever $i \neq j$, as the robot configuration at
time $t$. We formulate a {\em combined deploy and search} (CDS) strategy for multiple mobile robots as a modification to {\em sequential deploy and search} (SDS) strategy proposed in our previous work \cite{tase}. The fundamental concept is ``deploy" and ``search". In each search step, the uncertainty is reduced as
\[
\phi_{n+1}(q) = \phi_n(q)\min_i\{\beta(\| p_i - q
\|)\}\]
where, $n$ is the search step count, and $\beta: \mathbb{R}\mapsto [0,1]$ is the function acting as a
reduction factor for the uncertainty density $\phi$. The function
$\beta$ is such that, $1-\beta$ is the sensor effectiveness function.

We consider following multi-center objective function, which when maximized, maximizes a single step effectiveness
\begin{displaymath}
\mathcal{H}_n = \sum_{i\in \{1,2,\ldots,N\}}\int_{V_{i}}\phi_n(q)[1-\beta(\| p_i - q \|)]dQ
\end{displaymath}
where, $n$ is the search step count, $V_i$ is the Voronoi cell corresponding to the $i$-th agent/robot. The critical points of the objective function for a given $n$ is shown to be centroids of Voronoi cells with perceived density. Unlike SDS, the agents/robots do not wait till getting optimally deployed to perform search. The critical points are used only to get a direction of motion for robots. We show, that the CDS strategy can reduce the average uncertainty to any arbitrarily small value in finite time. The optimal deployment strategy has been analyzed in presence of some constraints on robot speed and limit on sensor range for convergence of the robot trajectories with the corresponding control laws responsible for the motion of robots.

Simulation experiments are carried out to validate and compare the performance of SDS and CDS
with two generic strategies, namely, {\em greedy} and {\em random} search. The simulation results indicate that
both the proposed search strategies perform quite well even when the conditions deviated from the assumed ones
such as sensor range limitations and the CDS strategy leads to somewhat shorter and
smoother trajectories than those of the SDS strategy with the same parameters.

\subsection{Organization of the papers} The paper is organized as follows. We preview a few mathematical concepts used in this work in Section 2. The multi-robot search problem is discussed in Section 3. Section 4 discusses the {\em combined deploy and search} (CDS) strategy. The objective function, instantaneous critical points, and control law responsible for robot motion are presented in this section. In Section 5 we impose a few constraints on robot speed and also limit on sensor range, propose control laws and provide convergence results. The CDS strategy is explained with help of illustrative examples in comparison with the SDS strategy in Section 6. A few implementation issues are discussed in Section 7. In section 8, we provide and discuss simulation results and the paper concludes in section 9.

\section{Mathematical preliminaries}In this section we preview
mathematical concepts from  Voronoi partition, LaSalle's invariance principle, and
Liebniz theorem used in the present work.

\subsection{Voronoi partition} Voronoi partition \cite{vor2,vor1}
is a widely used scheme of partitioning a given space and finds
applications in many fields such as CAD, image processing and sensor
coverage. We use the Voronoi decomposition scheme to partition the
search space. Here we briefly preview the concept.

By a \emph{partition} of a set $X$ we mean a collection of subsets
$W_i$ of $X$ with disjoint interiors such that their union is $X$
itself. Let $Q \subset \mathbb{R}^d$, be a convex polytope in $\mathbb{R}^d$, the $d$-dimensional Euclidean space. We
define the configuration space as $\mathbb{P} =\{\mathcal{P} =
\{p_1,p_2,\ldots,p_N\}\}$, $p_i \in Q$, position of $i$-th node in $Q$. The {\em Voronoi partition}
generated by $\mathcal{P} \in \mathbb{P}$ with respect to Euclidean
norm is the collection $\{V_i(\mathcal{P})\}_{i\in \{1,2,\dots,n\}}$
defined as,
\begin{displaymath}
V_i(\mathcal{P}) = \left \{ q \in Q | \parallel q-p_i \parallel \leq
\parallel q - p_j \parallel, \forall p_j \in \mathcal{P} \right\}
\end{displaymath}

The Voronoi cell $V_i$ is  the collection of those points which
are closest to $p_i$ compared to any other point in $\mathcal{P}$.
The boundary of each Voronoi cell is the union of a finite
number of line segments forming a closed $C^0$ curve. In
$\mathbb{R}^2$, The intersection of any two Voronoi cells is
either null, a line segment, or a point. In a general $d$-dimensional space, the
boundaries of the Voronoi cells are unions of convex subsets of at most $d-1$
dimensional hyperplanes in $\mathbb{R}^d$ and the intersection of
two Voronoi cells is either a convex subset of a hyperplane or a null set. Each of the
Voronoi cells is a topologically connected non-null set.

\subsection{LaSalle's Invariance principle} Here we state LaSalle's
invariance principle \cite{lasalle,lasalle2} used widely to study
the stability of nonlinear dynamical systems. We state the theorem
as in \cite{marquez} (Theorem 3.8 in \cite{marquez}).

Consider a dynamical system in a domain $D$
\begin{equation}
\label{dyn_lasalle} \dot x = f(x)\text{, } f:D \rightarrow
\mathbb{R}^d
\end{equation}

Let $V:D \rightarrow \mathbb{R}$ be a continuously differentiable
function and assume that
\begin{enumerate}
\item[(i)] $M \subset D$ is a compact set, invariant with respect to
the solutions of (\ref{dyn_lasalle}).

\item[(ii)] $\dot{V} \leq 0$ in $M$.

\item[(iii)] $E:\{x:x\in M\text{, and } \dot{V}(x) = 0 \}$; that is $E$
is set of all points of $M$ such that $\dot{V}(x)=0$.

\item[(iv)] $N$ is the largest invariant set in $E$.
\end{enumerate}

Then every solution of (\ref{dyn_lasalle}) starting in $M$
approaches $N$ as $ t \rightarrow \infty$.

Here by \emph{invariant set} we mean that if the trajectory is
within the set at some time, then it remains within the set for all
time. Important differences of the LaSalle's invariance principle as
compared to the Lyapunov Theorem are (i) $\dot{V}$ is required to be
negative semi-definite rather than negative definite and (ii) the
function $V$ need not be positive definite (see Remark on Theorem
3.8 in \cite{marquez}, pp 90-91).

\subsection{Leibniz Theorem and its generalization}
The Leibniz Theorem is widely used in fluid mechanics \cite{kundu},
and shows how to differentiate an integral whose integrand as well
as the limits of integration are functions of the variable with
respect to which differentiation is done. The theorem gives the
formula
\begin{equation}
\label{leibniz} \frac{d}{dy}\int_{a(y)}^{b(y)}F(x,y)dx = \int_a^b
\frac{\partial F}{\partial y}dx + \frac{db}{dy}F(b,y) -
\frac{da}{dy}F(a,y)
\end{equation}

Eqn. (\ref{leibniz}) can be generalized for a $d$-dimensional
Euclidean space as
\begin{equation}
\label{leibnizgen} \frac{d}{dy}\int_{\mathcal{V}(y)}F(x,y)d\mathcal{V}= \int_\mathcal{V}
\frac{\partial F}{\partial y}d\mathcal{V} + \int_{\mathbf{S}}
\mathbf{n}(x).\mathbf{u}(x)FdS
\end{equation}
where, $\mathcal{V} \subset \mathbb{R}^d$ is  the volume in which the
integration is carried out, $d\mathcal{V}$ is the differential
volume element, $\mathbf{S}$ is the bounding hypersurface of $\mathcal{V}$,
$\mathbf{n}(x)$ is the unit outward normal to $\mathbf{S}$ and
$\mathbf{u}(x) = \frac{d\mathbf{S}}{dy}(x)$ is the rate at which the surface
moves with respect to $y$ at $x \in \mathbf{S}$.

\section{Multi-robot search} In this section we discuss the problem
addressed in this paper. $N$ robots perform search operation in an
unknown environment. The lack of information is modeled as an
uncertainty density distribution over the search space $Q$. The
problem addressed in this paper is that of deploying $N$ robots in
$Q$ to collect information, thereby reducing the uncertainty density
distribution over $Q$. The problem formulation is stated formally as

\begin{enumerate}
\item $Q \subset \mathbb{R}^d$ is a convex polytope and is the search
space.

\item $\phi:Q\mapsto[0,1]$ defines the uncertainty density function
representing lack of information.

\item $N$ robots, equipped with sensors and communication
equipment, deploy themselves in $Q$, and gather information, thus
reducing the uncertainty.

\item $P(t) = \{p_1(t),p_2(t),\ldots,p_N(t)\} \subset Q$, with $p_i \neq p_j$ whenever $i \neq j$, denotes the configuration
of the multi-robot system at time $t$, $p_i(t)$ denotes the position
of the $i$-th robot at time $t$. In future, for convenience, we drop
the variable $t$ and refer to the positions by just $p_i$.

\item Sensor's effectiveness at a point reduces with distance from the sensor.

\item Ideally, we are looking for an optimal way of utilizing the robots to acquire
complete information about $Q$, and thus have $\phi(q) = 0$,
$\forall q \in Q$.
\end{enumerate}

During the search operation, sensors gather information about $Q$, reducing the uncertainty density as,
\begin{equation}
\label{phiupdate} \phi_{n+1}(q) = \phi_n(q)\min_i\{\beta(\| p_i - q
\|)\}
\end{equation}
where, $\phi_n(q)$ is the density function at the $n$-th iteration;
$\beta: \mathbb{R}\mapsto [0,1]$ is a function of Euclidian distance
of a given point in space from the robot, and acts as the factor of
reduction in uncertainty by the sensors; and $p_i$ is the position
of the $i$-th sensor. At a given $q \in Q$, only the robot with the
smallest $\beta(\| p_i - q \|)$, that is, the robot which can reduce
the uncertainty by the largest amount, is active. If any robot searches within its Voronoi cell, then the updating function (\ref{phiupdate}) gets implemented automatically, That is, the function $\min_i\{\beta(\|p_i - q \|)\}$ is simply $\beta(\|p_i - q \|)$, where $p_i \in V_i$.

In the SDS strategy proposed in \cite{tase}, the agents get optimally deployed before performing search. In order to maximize the search effectiveness in each search step, following objective function was considered to be maximized.
\begin{equation}
\label{obj1}
\begin{array}{lcl}
\mathcal{H}_n &=& \int_Q \Delta\phi_n(q)dQ\\
              &=& \int_Q \max_{i\in\{1,2,\ldots,N\}}\{(|\phi_n(q) - \beta(\|
p_i - q \|)\phi_n(q)|)\}dQ\\
              &=& \int_Q (\phi_n(q) - \min_{i\in\{1,2,\ldots,N\}}\{\beta(\|
p_i - q \|)\}\phi_n(q))dQ\\
              &=& \sum_{i\in\{1,2,\ldots,N\}}\int_{V_{i}}\phi_n(q)(1-\beta(\|
p_i - q \|))dQ
\end{array}
\end{equation}
where, $n$ is ``deploy" and ``search" count, $V_{i}$ is the Voronoi cell corresponding to the $i$-th agent, and $p_{i} \in Q$ is the position
of the $i$-th agent.

The search effectiveness function $\beta: \mathbb{R}\mapsto [0,1]$
is a non-decreasing function capturing effectiveness of the sensor.
%A continuously differentiable function leads to a possible
%closed-form solution to the optimization problem (\ref{obj1}).
Consider
\begin{displaymath}
\beta(r) = 1- ke^{-\alpha r^2},\quad k \in (0, 1)\quad
\text{and $\alpha > 0$}
\end{displaymath}

\begin{figure}[h]
\centerline{\psfig{figure=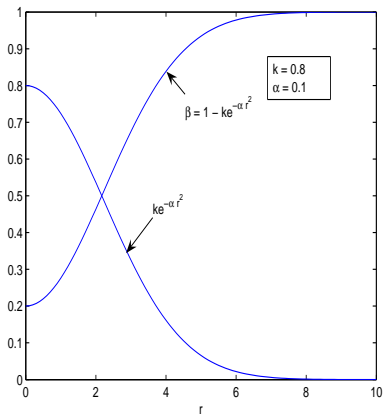,height=6cm,width=6cm}}
\caption{\small The sensor effectiveness and updating
function}\label{beta}
\end{figure}

Here, $ke^{-\alpha r^2}$ represents the effectiveness of the sensor
which is maximum at $r=0$ and tends to zero as $r \rightarrow
\infty$ and $\beta$ is minimum at $r = 0$ (effecting maximum
reduction in $\phi$) and tends to unity as $r \rightarrow \infty$
(Figure \ref{beta}) (change in $\phi$ reduces to zero as $r$
increases). Most sensors' effectiveness reduces over distance as the signal to noise ratio increases with the distance. Thus $\beta$, which is upside down Gaussian, can model a wide variety of sensors with two tunable parameters $k$ and $\alpha$.

The optimal deployment configuration was shown to be a variation of centroidal Voronoi configuration, where each robot is located at the centroid of its Voronoi cell computed with a density $\tilde{\phi}_n(q) =
\phi_n(q)ke^{-\alpha r_i^2}$, which is the density as perceived by the
sensor, with $\beta(r) = 1- ke^{-\alpha r^2}$, $k \in (0, 1)$, and $\alpha > 0$, as search effectiveness function.

We use Voronoi partition in formulating the search strategies which along with advantages can cause
some computation overhead (see \cite{franz} and references therein for time complexity of computing Voronoi partitions). This issue has been addressed in the literature (see \cite{bullo1} and references therein) and there
are a few algorithms that efficiently implement Voronoi partition related computations. Also, Voronoi based strategies result
in collision free trajectories in a natural way, which is an added advantage.

\subsection{Optimal deployment strategy} We have shown that the uncertainty reduction will be maximized in a single step of search, if the agents are located at the centroids of respective Voronoi cells. In SDS, the agents get deployed optimally in this sense before performing search. In this section we will provide the control law to make the agent move toward the centroids and achieve the optimal deployment configuration. Though we have discussed the control law in \cite{tase}, we provide it here for the sake of completeness and clarity.

Typically search problems do not consider dynamics of search agents as the focus is more on the effectiveness
of search, that is, being able to identify region of high uncertainty and distribute search effort to reduce uncertainty.
Moreover, it is usually assumed that the search region is large compared to the physical size of the agent or the
area needed for the agent to maneuver. In this paper, we assume that the agents are modeled as simple first dynamical systems as
\begin{equation}
\label{dyn1} \dot p_{i} = u_i
\end{equation}

Consider the control law
\begin{equation}
\label{ctrl1} u_i = -k_{prop}(p_{i} - \tilde{C}_{V_{i}})
\end{equation}

Control law (\ref{ctrl1}) makes the agents move toward
$\tilde{C}_{V_{i}}$ for positive control gain, $k_{prop}$.

We have shown in \cite{tase}, using LaSalle's invariance principle, that the trajectories of the agents governed by the
control law (\ref{ctrl1}), starting from any initial condition
$\mathcal{P}(0) \in Q^N$, will asymptotically converge to the
critical points of $\mathcal{H}_n$.

\section{Combined deploy and search (CDS) strategy}
In the  SDS strategy, The robots first get optimally deployed and then search is performed. The ``deploy" and ``search" steps continue till the uncertainty density is reduced below a desired value. Here the optimal deployment strategy ensures that the uncertainty density reduction is maximized in each search step. But it does not guarantee optimal trajectories of the robot. During the deployment stage, the robots move without utilizing the sensors. Intuitively, it seems that the trajectories
will be closer to optimal if, as the robots are moving toward the
respective $\tilde{C}_{V_i}$, they also
simultaneously perform the search operation in discrete steps. We define the {\em latency}, $t_s$, of the
robots as the maximum time taken to acquire the information, process
it, and successfully update the uncertainty density. The time interval between each search should be more than $t_s$. Here we shall formulate
such a strategy and name it {\em combined deploy and search} (CDS) strategy.

\subsection{Density update} Here we provide the problem formulation for
the {\em combined deploy and search} strategy. Assume that the
index n represents the intermediate step at which the search is performed and uncertainty
density is updated. Using the uncertainty density update rule (\ref{phiupdate}) discussed earlier we can get,
\begin{equation}
\label{ContSrchDiffrence} \Delta_n \phi(q) = \phi_{n+1}(q) -
\phi_n(q) = \phi_n(q) \min_i (1 - \beta(\parallel p_i - q
\parallel))
\end{equation}
\begin{equation}
\label{bigphi} \Phi_n = \int_Q \phi_n(q)dQ
\end{equation}
Integrating (\ref{ContSrchDiffrence}) over $Q$,
\begin{equation}
\label{contsrch2} \Delta \Phi_n = \sum_{i\in\{1,2,\ldots,N\}}
\int_{V_i}\phi_n(q)(1 - \beta(\parallel p_i - q
\parallel))dQ
\end{equation}

\subsection{Objective function}
The objective function (\ref{obj1}), used for SDS strategy \cite{tase}, is fixed for each deployment step as $\phi_n(q)$ is
fixed for the $n$-th iteration. In {\em combined deploy and search},
the search task is performed as the robots move. Now an objective function to be maximized in order to maximize the uncertainty reduction at the $n$-th search step is
\begin{equation}
\label{obj_cs1} \mathcal{H}_n = \Delta \Phi_n
               = \sum_{i\in\{1,2,\ldots,N\}}
\int_{V_i}\phi_n(q)(1 - \beta(\parallel p_i - q
\parallel))dQ
\end{equation}
Note that the above objective function is same as (\ref{obj1}) except for the fact that $n$ in this case represents the search step count, whereas in (\ref{obj1}) it represents `deploy and search' step count.
For $\beta(r) = 1- ke^{-\alpha r^2}$, the objective function
(\ref{obj_cs1}) becomes,
\begin{equation}
\label{obj_cs2} \mathcal{H}_n
=\sum_{i\in\{1,2,\ldots,N\}}
\int_{V_i}\phi_n(q)ke^{-\alpha r_i^2}dQ
\end{equation}

It can be noted that for a given $n$, the uncertainty density $\phi_n(q)$ at any $q \in Q$ is constant. The gradient is given as (see Theorem A.1 in Appendix),
\begin{eqnarray}
\label{grad_cs}
\frac{\partial\mathcal{H}_n}{\partial p_i} &=&
\sum_{i\in\{1,2,\ldots,N\}}\int_{V_{i}}\phi_n(q)ke^{-\alpha (\| p_i - q\|)^2}(-2\alpha)(p_{i}-q)dQ \nn \\
&=& -2\alpha \tilde{M}_{V_{i}}(p_{i} - \tilde{C}_{V_{i}})
\end{eqnarray}
where $\tilde{M}_{V_{i}}$ and $\tilde{C}_{V_{i}}$ are the mass and
the centroid of $V_{i}$ with respect to $\tilde{\phi}_n(q) =
\phi_n(q)ke^{-\alpha r_i^2}$, which is the density as perceived by the
sensor. The critical points are same as those obtained for SDS. But the uncertainty changes in every time step and hence the critical points also change. Hence, the corresponding critical points are only the instantaneous critical points. It should be noted that the above treatment is valid for any non-decreasing continuously differentiable $\beta(\cdot)$ with $\tilde{\phi}(\cdot)$ depending on exact nature of the function $\beta(\cdot)$.
We use the control law
\begin{equation}
\label{ctrl2}u_i = -k_{prop}(p_{i} - \tilde{C}_{V_{i}})
\end{equation}
Control law (\ref{ctrl1}) makes the robots move toward
$\tilde{C}_{V_{i}}$ for positive control gain, $k_{prop}$.

The instantaneous critical points and the gradient (\ref{grad_cs}) are used in control law (\ref{ctrl2}) only to make the robots move toward the instantaneous centroids rather than deploying them optimally. Thus, it is not possible to prove any optimality of deployment and we do not prove the convergence of the trajectories here in case of CDS. In CDS, compared to SDS, robots perform more frequent searches instead of waiting till the optimal deployment maximizing per step uncertainty reduction.

 To implement the control law, centroid of each Voronoi cell needs to be computed.  The computational overhead of computing the centroid can be reduced at the cost of slower convergence using methods reported in the literature such as random sampling and stochastic approximation \cite{lvq,pages}. In addition, we discretize the search space into grids while implementing the strategy. This simplifies the computation of the centroid of Voronoi cells. The main focus of this paper is design and demonstration of the multi-robot search strategy and finer issues such as computation complexities are beyond the scope of this paper.

It can be shown that he {\em combined deploy
and search} strategy is spatially distributed over the Delaunay
graph $\mathcal{G}_D$. Here by spatially distributed we mean that information from neighboring robots is sufficient for computation of control input. A Delaunay graph $\mathcal{G}_D$ is an undirected graph, where two agents/robots are said to be neighbors (connected by an edge) if the corresponding Voronoi cells are adjacent.

\begin{thm}
\label{CS_convergence} The {\em combined deploy
and search} strategy can reduce the average uncertainty to any
arbitrarily small value in finite time.
\end{thm}

\noindent {\em Proof.~} Consider the  uncertainty density update law
(\ref{phiupdate}) for any $q \in Q$,
\begin{equation}
\label{phiupdate_n-1} \phi_n(q) = (1 - ke^{-\alpha
{r_i}^2})\phi_{n-1}(q)= \gamma_{n-1}\phi_{n-1}(q)
\end{equation}
where, $r_i$ is the distance of point $q \in Q$ from the $i$-th
robot, such that $q \in V_i$, the Voronoi cell corresponding to
it and, $\gamma_{n-1} = (1 - ke^{-\alpha {r_i}^2})$.

Applying the above update rule recursively, we have,
\begin{equation}
\label{phiupdate_recur} \phi_n(q) = \gamma_{n-1}\gamma_{n-2}\ldots
\gamma_1\gamma_0\phi_0(q)
\end{equation}

Let $D(Q) := \max_{p,q \in Q}(\parallel p - q \parallel)$. We note
that
\begin{enumerate}
\item[(i)] $0 < k < 1$
\item[(ii)] $0 \leq r_i \leq D(Q)$. $D(Q)$ is bounded as the set $Q$
is bounded.
\item[(iii)] $0 \leq \gamma_j \leq 1 - ke^{-\alpha \{D(Q)^2\}} = l$ (say), $j \in \mathbb{N}$;
and $l < 1$
\end{enumerate}

Now consider the sequence $\{\Gamma\}$ ,
\begin{displaymath}
\Gamma_n = \gamma_n\gamma_{n-1}\ldots \gamma_1\gamma_0 \leq l^{n+1}
\end{displaymath}

Taking limits as $n \rightarrow \infty$,
\begin{displaymath}
\lim_{n \rightarrow \infty}\Gamma_n \leq \lim_{n \rightarrow
\infty}l^{n+1} = 0
\end{displaymath}

Thus,
\begin{displaymath}
\lim_{n \rightarrow \infty}\phi_{n}(q) = \lim_{n \rightarrow
\infty}\Gamma_{n-1}\phi_0(q) = 0
\end{displaymath}
As the uncertainty density $\phi$ vanishes at each point $q \in Q$
in the limit, the average uncertainty density over $Q$ is bound to
vanish as $n \rightarrow \infty$. Thus, the average uncertainty
density can be reduced to arbitrarily small
value in finite time. \hfill $\Box$\\

It can be observed that the above proof does not depend on the control law. The theorem depends only on the
outcome of the choice of the updating function (\ref{phiupdate}), along with the fact that there is no sensor range limitation, and
that the search space Q is bounded. In addition, the theorem does not address the issue of optimality of the strategy
which, in fact, depends on the control law which is responsible for the motion of the robots. Further, unlike SDS, maximal uncertainty reduction is also not guaranteed in each search step. In case of SDS, the reduction in the uncertainty in each step in SDS is
\begin{equation}
\label{SDSrate}
\mathcal{H}_n^* = \sum_i \int_{V_i} \phi_n(q)ke^{-\alpha(\|\tilde{C}_{V_i} - q\||)^2} dQ
\end{equation}
which is the maximum possible reduction in a single step. The deployment is such that uncertainty will be reduced
to a maximum possible extent in a step, given by the above formula. whereas in CDS the uncertainty reduction in $n$-th search step is given by
\begin{equation}
\label{SDSrate1}
\mathcal{H}_n^* = \sum_i \int_{V_i} \phi_n(q)ke^{-\alpha(\|p_i - q\||)^2} dQ
\end{equation}
where it is not required that $p_i = \tilde{C}_{V_i}$ while performing search. Though the uncertainty reduction in a given search step $n$ in CDS is less than that in SDS, ss will be seen in later sections, the CDS performs better compared to SDS in terms of faster uncertainty reduction due to more frequent searches.

\section{Constraints on robot speed and sensor range}In previous sections we have formulated the multi-robot search strategies under ideal conditions and provided a few useful analytical results for convergence and spatial distributedness providing an analytically sound platform for
analysis. But in a practical situation these conditions may be violated. It is more likely that the robots will have
limit on maximum speed or they may be constrained to move with a constant speed. Further, the sensors could
have limit on their range. In this section we analyze the proposed strategies in the presence of speed and sensor
range limitations.

\subsection{Maximum speed constraint} Let the robots have a constraint
on maximum speed of ${U_{max}}_i$, for $i=1,\ldots,n$. Now consider
the control law

\begin{equation}
\label{ctrl_with_sat} u_i =
\begin{cases}
-k_{prop}(p_{i} - \tilde{C}_{V_{i}}) & \text{If $u_i \leq {U_{max}}_i$} \\
-{U_{max}}_i\frac{(p_{i} - \tilde{C}_{V_{i}})}{\parallel(p_{i} -
\tilde{C}_{V_{i}})\parallel} & \text{Otherwise}
\end{cases}
\end{equation}

The control law (\ref{ctrl_with_sat}) makes the robots move toward
their respective centroids with saturation on speed.

\begin{thm}
\label{saturation_speed_stability} The trajectories of the robots
governed by the control law (\ref{ctrl_with_sat}), starting from any
initial condition $\mathcal{P}(0) \in Q^N$, will asymptotically
converge to the critical points of $\mathcal{H}_n$.
\end{thm}

\noindent {\it Proof.~} Consider  $V(\mathcal{P}) = -\mathcal{H}_n$.

\begin{equation}
\label{vdot_sat}
\begin{array}{lcl}
{\dot V}(\mathcal{P}) &=& -\frac{d\mathcal{H}_n}{dt} = -\sum_{i\in\{1,2,\ldots,N\}} \frac{\partial \mathcal{H}_n}{\partial p_{i}}\dot{p}_{i}\\
&=& \begin{cases} 2\alpha \sum_{i\in\{1,2,\ldots,N\}} \tilde{M}_{V_{I}}(p_{i} -
 \tilde{C}_{V_{i}})(-k_{prop})(p_{i} - \tilde{C}_{V_{i}}) \text{, ~~If $u_i \leq {U_{max}}_i$} \\
2\alpha \sum_{i\in\{1,2,\ldots,N\}} \tilde{M}_{V_{I}}(p_{i} -
 \tilde{C}_{V_{i}})(-{U_{max}}_i)\frac{(p_{i} -
C'_{V_{i}})}{(\|p_{i} - \tilde{C}_{V_{i}}\|)}\text{, ~~otherwise}\\
\end{cases} \\
&=&
\begin{cases}
-2\alpha k_{prop}\sum_{i\in\{1,2,\ldots,N\}}  \tilde{M}_{V_{i}}(\|p_{i} -
\tilde{C}_{V_{i}}\|)^2  \text{, ~~If $u_i \leq {U_{max}}_i$}  \\
-2\alpha \sum_{i\in\{1,2,\ldots,N\}} {U_{max}}_i\tilde{M}_{V_{i}}\frac{(\|p_{i} -
\tilde{C}_{V_{i}}\|)^2}{(\parallel p_{i} - \tilde{C}_{V_{i}}\parallel)}
 \text{, ~~ otherwise}
\end{cases}
\end{array}
\end{equation}

We observe that
\begin{enumerate}
\item $V: Q\mapsto \mathbb{R}$ is continuously differentiable in
$Q$ as $\{V_i\}$ depends at least continuously on $\mathcal{P}$, and $\dot V$ is continuous as $u$ is continuous if not smooth.

\item $M = Q$ is a compact invariant set.

\item ${\dot V}$ is negative definite in $M$.

\item $E = \dot{V}^{-1}(0) = \{\tilde{C}_{V_{i}}\}$.

\item $E$ itself is the largest invariant subset of $E$ by the control law
(\ref{ctrl_with_sat}).

\end{enumerate}

Thus, by LaSalle's invariance principle, the trajectories of the
robots governed by control law (\ref{ctrl_with_sat}), starting from
any initial configuration $\mathcal{P}(0) \in Q^N$, will
asymptotically converge to the set $E$, the critical points of
$\mathcal{H}_n$, that is, the centroidal Voronoi configuration with
respect to the density function as perceived by the sensors. \hfill
$\Box$

\subsection{Constant speed control} The robots may have a constraint
of moving with a constant speed $U_i$. But we let the robots slow down as they approach the critical points. For $i=1,\ldots,n$, consider
the control law
\begin{equation}
\label{const_speed_ctrl} u_i =
\begin{cases} -U_i\frac{(p_{i} -
\tilde{C}_{V_{i}})}{\parallel(p_{i} - \tilde{C}_{V_{i}})\parallel}
\text{, if $\|p_i - \tilde{C}_{V_i}\| \geq \delta$} \\
-U_i(p_i - \tilde{C}_{V_i})/\delta \text{, otherwise} \end{cases}
\end{equation}
where, $\delta > 0$, predefined value, such that the control law (\ref{const_speed_ctrl}) makes the robots move
toward their respective centroids with a constant speed of $U_i$ when they at a distance greater than $\delta$ from the corresponding centroids and slow down as they approach them.

\begin{thm}
\label{const_speed_stability} The trajectories of the robots
governed by the control law (\ref{const_speed_ctrl}), starting from
any initial condition $\mathcal{P}(0) \in Q^N$, will asymptotically
converge to the critical points of $\mathcal{H}_n$.
\end{thm}

\noindent {\it Proof}. Consider  $V(\mathcal{P}) = -\mathcal{H}_n$.

\begin{equation}
\label{vdot_const}
{\dot V}(\mathcal{P}) = \begin{cases} -2\alpha \sum_{i\in\{1,2,\ldots,N\}}U_i
\tilde{M}_{V_{i}}\frac{(\|p_{i} -
\tilde{C}_{V_{i}}\|)^2}{(\parallel p_{i} -
\tilde{C}_{V_{i}}\parallel)},\\
 \text{~~~~~ if $\|p_i - \tilde{C}_{V_i}\| \geq \delta$} \\
-2\alpha \sum_{i\in\{1,2,\ldots,N\}}U_i \tilde{M}_{V_{i}}(p_i - \tilde{C}_{V_i})/\delta,\\ \text{~~~~~otherwise}
\end{cases}
\end{equation}

We observe that
\begin{enumerate}
\item $V: Q\mapsto \mathbb{R}$ is continuously differentiable in
$Q$ as $\{V_i\}$ depends at least continuously on $\mathcal{P}$, and $\dot V$ is continuous as $u$ is continuous if not smooth.

\item $M = Q$ is a compact invariant set.

\item ${\dot V}$ is negative definite in $M$.

\item $E = \dot{V}^{-1}(0) = \{\tilde{C}_{V_{i}}\}$.

\item $E$ itself is the largest invariant subset of $E$ by the control law
(\ref{const_speed_ctrl}).

\end{enumerate}

Thus, by LaSalle's invariance principle, the trajectories of the
robots governed by control law (\ref{const_speed_ctrl}), starting
from any initial configuration $\mathcal{P}(0) \in Q^N$, will
asymptotically converge to the set $E$, the critical points of
$\mathcal{H}_n$, that is, the centroidal Voronoi configuration with
respect to the density function as perceived by the sensors. \hfill
$\Box$

\subsection{Effect of sensor range limits}

In reality, it is likely that the sensors will have limitations on
their range. The sensors, in addition to having a monotonically
decreasing effectiveness with Euclidean distance, might be totally
insensitive to signals at distances larger than $R$, the sensor
range limit. It can also be thought of as follows: when the
effectiveness falls below, say, a certain threshold value, for all
practical purposes, it can be assumed to be ineffective.

In order to incorporate the sensor range limit, we need to modify
the objective function (\ref{obj1}). We do so by modifying $\beta$
suitably.

\begin{prop}
\label{saturation1} Let $\beta(r)$ and $\hat\beta(r) = c +
\beta(r)$, where $c$ is a real constant, be two sensor detection
functions. The corresponding objective functions
$\mathcal{H}(\mathcal{P})$ and $\hat{\mathcal{H}}(\mathcal{P})$,
respectively, have the same critical points.
\end{prop}

\noindent {\it Proof.} The objective function with sensor detection
function $\hat \beta$ is
\begin{displaymath}
\begin{array}{lcl}
\hat{\mathcal{H}}(\mathcal{P}) &=&
\sum_{i\in\{1,2,\ldots,N\}}\int_{V_{i}}\phi_n(q)(1-\hat \beta(\| p_i -
q \|))dQ\\
&=& \sum_{i\in\{1,2,\ldots,N\}}\int_{V_{i}}\phi_n(q)(1-\beta(\| p_i - q \|) -c)dQ\\
&=& \sum_{i\in\{1,2,\ldots,N\}}\int_{V_{i}}\phi_n(q)(1-\beta(\| p_i - q \|))dQ
-c\sum_{i\in\{1,2,\ldots,N\}}\int_{V_{i}}\phi_n(q)dQ\\
&=& \mathcal{H} - c \sum_{i\in\{1,2,\ldots,N\}}\int_{V_{i}}\phi_n(q)dQ
\end{array}
\end{displaymath}
The second term not being a function of $p_i$, we have,
\begin{displaymath}
\frac{\partial(\hat{\mathcal{H}})}{\partial p_i}(\mathcal{P}) =
\frac{\partial(\mathcal{H})}{\partial p_i}(\mathcal{P})
\end{displaymath}
\mbox{}\hfill $\Box$

\begin{figure}
\centerline{\psfig{figure=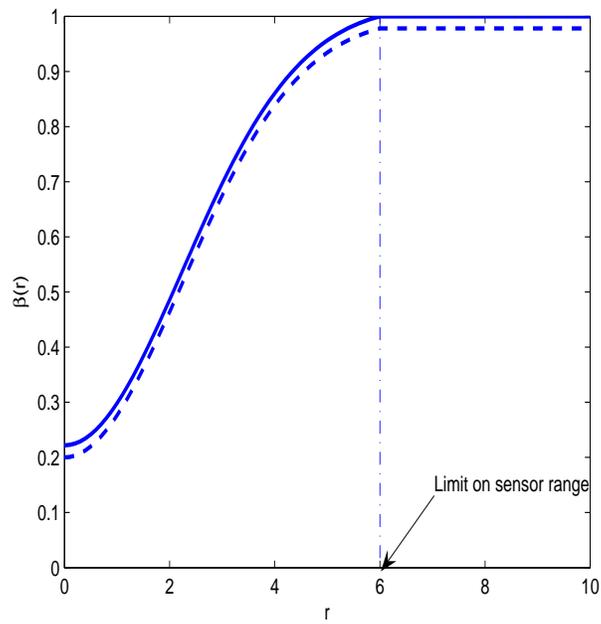,height=9cm,width=9cm}}
\caption{Illustration of $\tilde{\beta}$ and $\hat{\beta}$ in presence of the limit on sensor range. The dotted curve represent the sensor detection function $\tilde{\beta}$ and solid curve is the actual sensor detection function $\hat{\beta}(r) = \tilde{\beta}(r) + (1- \beta(R))$}
\label{sensor_range_fig}
\end{figure}

Consider the objective function with a saturation on $\beta$. Let
\begin{equation}
\label{beta_sat} \tilde{\beta}(r) = \begin{cases}\beta(r),
\text{if $r < R$} \\
\beta(R),  \text{otherwise}
\end{cases}
\end{equation}
The sensor detection function is shown in Figure \ref{sensor_range_fig} with dotted curve for $R=6$. For $r > 6$, the sensor detection function remains fixed at $\beta(6)$.

Consider the objective function
\begin{equation}
\label{obj_sensr_rng} \tilde{\mathcal{H}} =
\sum_{i\in\{1,2,\ldots,N\}}\int_{(V_{i}\cap\bar{B}(p_i,R))}\phi_n(q)(1-\tilde{\beta}(\|
p_i - q \|))dQ
\end{equation}

It is easy to show that the gradient of the objective function with
the new update function $\tilde{\beta}$  is given by,

\begin{equation}
\label{grad_sensr_rng} \frac{\partial(\tilde{\mathcal{H}})}{\partial
p_i}(\mathcal{P}) =
2\tilde{M}_{(V_{i}\cap\bar{B}(p_i,R))}(\tilde{C}_{(V_{i}\cap\bar{B}(p_i,R))}
- p_i)
\end{equation}
where, the mass $\tilde{M}$ and the centroid $\tilde{C}$ are now
computed within the region $(V_{i}\cap\bar{B}(p_i,R))$, that is, the
region of Voronoi cell $V_i$, which is accessible to the $i$-th
robot. The critical points are  $p_i =
\tilde{C}_{(V_{i}\cap\bar{B}(p_i,R))}$.

However, in reality, the detection function should become one
beyond the range $R$ (that is there is no reduction in uncertainty at points which are at a distance more than $R$ from the sensor. To achieve this, let
\begin{eqnarray}
 \hat \beta(r) &=&
\begin{cases}
\beta(r) + (1-\beta(R)) = 1-k(e^{-\alpha r^2} - e^{-\alpha R^2}),\text{~if~} r<R \\
1, ~\text{otherwise}
\end{cases}\\
& = & \tilde \beta(r) + (1 - \beta(R)) \nn
\end{eqnarray}
Figure \ref{sensor_range_fig} shows the sensor detection function $\hat{\beta}$ in solid line for $R=6$. After $r=6$, the limit on sensor range, the sensor is ineffective indicated by the value 1. There will not be any reduction in uncertainty density in this region $r > 6$.

Now by Proposition \ref{saturation1},
$\hat{\mathcal{H}}(\mathcal{P})$, the objective function
corresponding to the detection function $\hat \beta$, has the same
critical points as those of the objective function
$\tilde{\mathcal{H}}(\mathcal{P})$.

The control law making robots move toward the new critical point is

\begin{equation}
\label{ctrl_sensr_rng} u_i = -k_{prop}(p_{i} -
\tilde{C}_{(V_{i}\cap\bar{B}(p_i,R))})
\end{equation}

\noindent {\it Remark 4.~} The control law (\ref{ctrl_sensr_rng}) is
spatially distributed under the $r$-limited Delaunay graph
$\mathcal{G}_{LD}$, for any robot configuration $\mathcal{P}$.

\begin{thm}
\label{ctrl_senr_rng_stability} The trajectories of the robots
governed by the control law (\ref{ctrl_sensr_rng}), starting from
any initial condition $\mathcal{P}(0) \in Q^N$, will asymptotically
converge to the critical points of $\hat{\mathcal{H}}$.
\end{thm}

\noindent {\it Proof.~} The proof is based on LaSalle's Theorem similar to those of Theorems \ref{saturation_speed_stability} and \ref{const_speed_stability} with $V = - \tilde{\mathcal{H}}(\mathcal{P})$. It can be shown that $V$ is continuously differentiable using Theorem 2.2 in \cite{bullo2}. \hfill $\Box$

\section{Illustrative examples}In this section we show some
simulation results to illustrate the CDS strategy in comparison with SDS \cite{tase}. More detailed simulation results will be
presented in a later section.

\begin{figure}
\centerline{
\subfigure[]{\psfig{figure=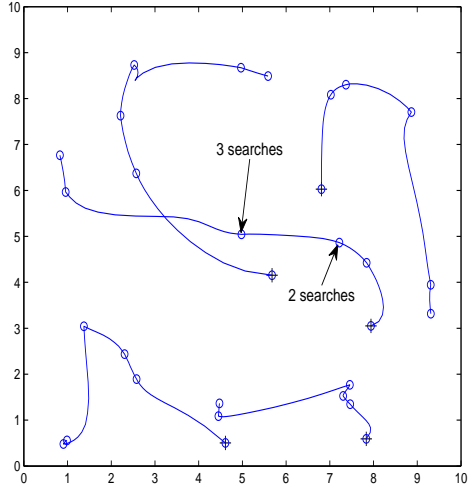,height=7.5cm,width=7.5cm}\label{dns2}}
\subfigure[]{\psfig{figure=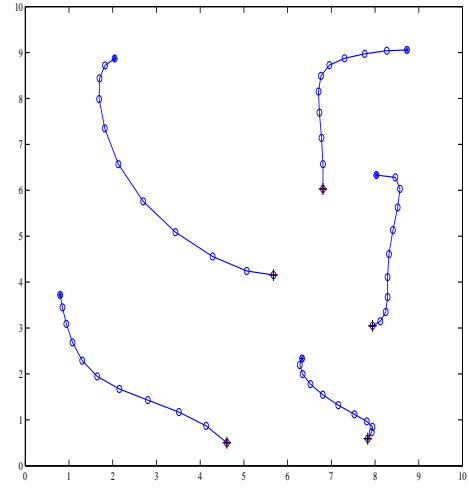,height=7.5cm,width=7.5cm}\label{cs2}}}
\centerline{
\subfigure[]{\psfig{figure=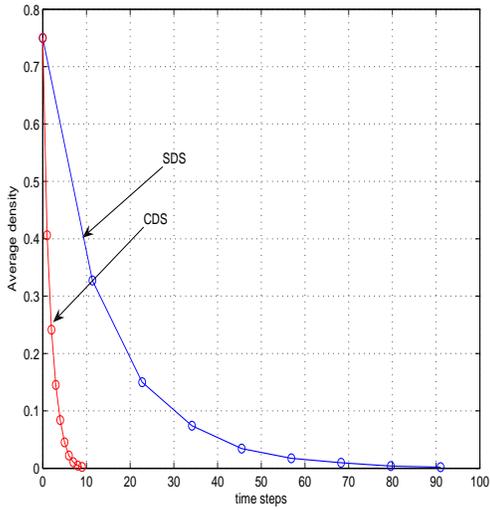,height=7.5cm,width=7.5cm}}
\subfigure[]{\psfig{figure=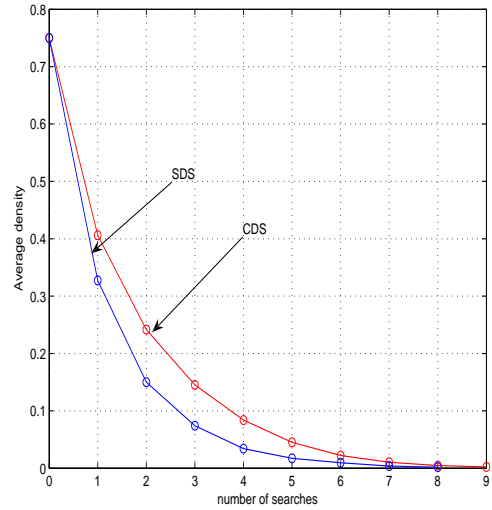,height=7.5cm,width=7.5cm}}
} \caption{Trajectories of robots with $N$=5 and without limit on
the sensor range for (a) SDS strategy and (b) CDS strategy. In
both cases the points marked `+' indicate the starting locations
of robots and `o' indicate the end of deployment and points in
space where search is being performed. In SDS, at a few places the
search is performed more than once. This is indicated in
trajectory of `robot 1'. The reduction in average uncertainty
density is shown in (c) against the number of time steps and (d)
against the number of searches, for SDS and CDS. Even in (c) and
(d), `o' indicate the search instances.}\label{comp2}
\end{figure}

Figure \ref{comp2} (a) and (b) compares the trajectories of robots
with SDS and CDS strategies. The trajectories with CDS are much
smoother and shorter. The instances of search are indicated by `o'
along the trajectories. It can be seen that the search is
performed at every discrete step in CDS, whereas the search is
performed only after each optimal deployment SDS.
Though there are 8 ``deploy and search" steps in SDS, only 5 `o's
are visible. In two of steps, multiple searches have been
performed as the centroids in successive steps were closer than some tolerance limit
$d_{tol} = 0.3$. Thus, there was no movement in corresponding
deployment step.

Figure \ref{comp2}(c) compares the history of uncertainty density
of SDS and CDS, and it can be observed that the CDS reduces
uncertainty relatively faster than SDS, in terms of number of time
steps. Figure \ref{comp2} (d) shows the reduction in average
uncertainty density with number of searches for SDS and CDS. It
can be observed from this figure that SDS reduces the uncertainty
in relatively fewer steps. This is apparent by very concept of
optimal deployment in SDS. CDS takes about 4 searches to reduce
uncertainty below 0.1, whereas SDS does this in only 3 searches.
If SDS requires over 30 time steps to achieve this reduction, CDS
needs 4 time steps. We can observe a tradeoff between the number
of searches and and number of time steps required to accomplish in
CDS and SDS. Once the uncertainty
reduces to a large extent in initial search steps, by nature of
the uncertainty density update rule (\ref{phiupdate}), amount of
reduction in subsequent searches is less in both SDS and CDS.

\begin{figure}
\centerline{
\subfigure[]{
\psfig{figure=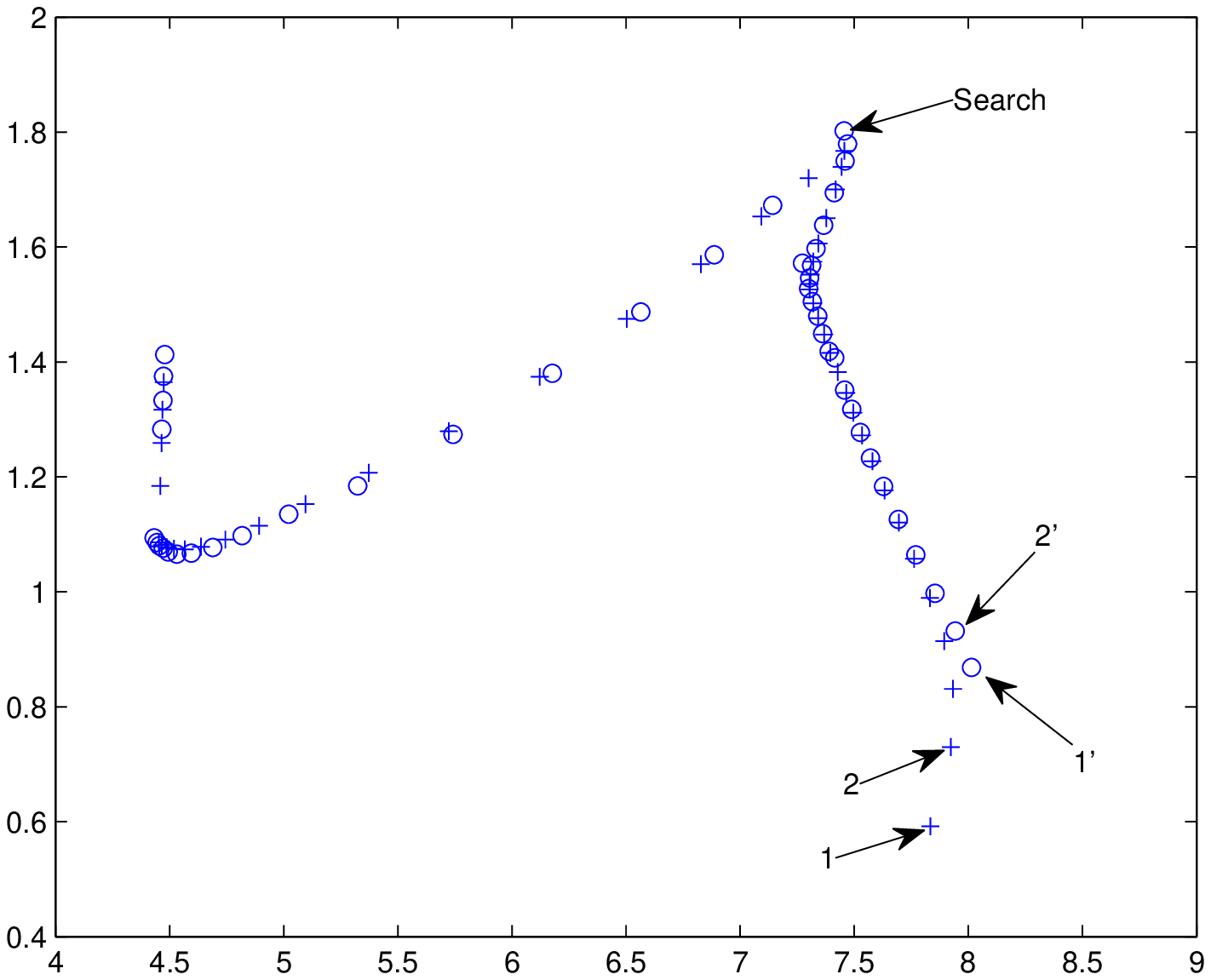,height=7.5cm,width=7.5cm}}
\subfigure[]{\psfig{figure=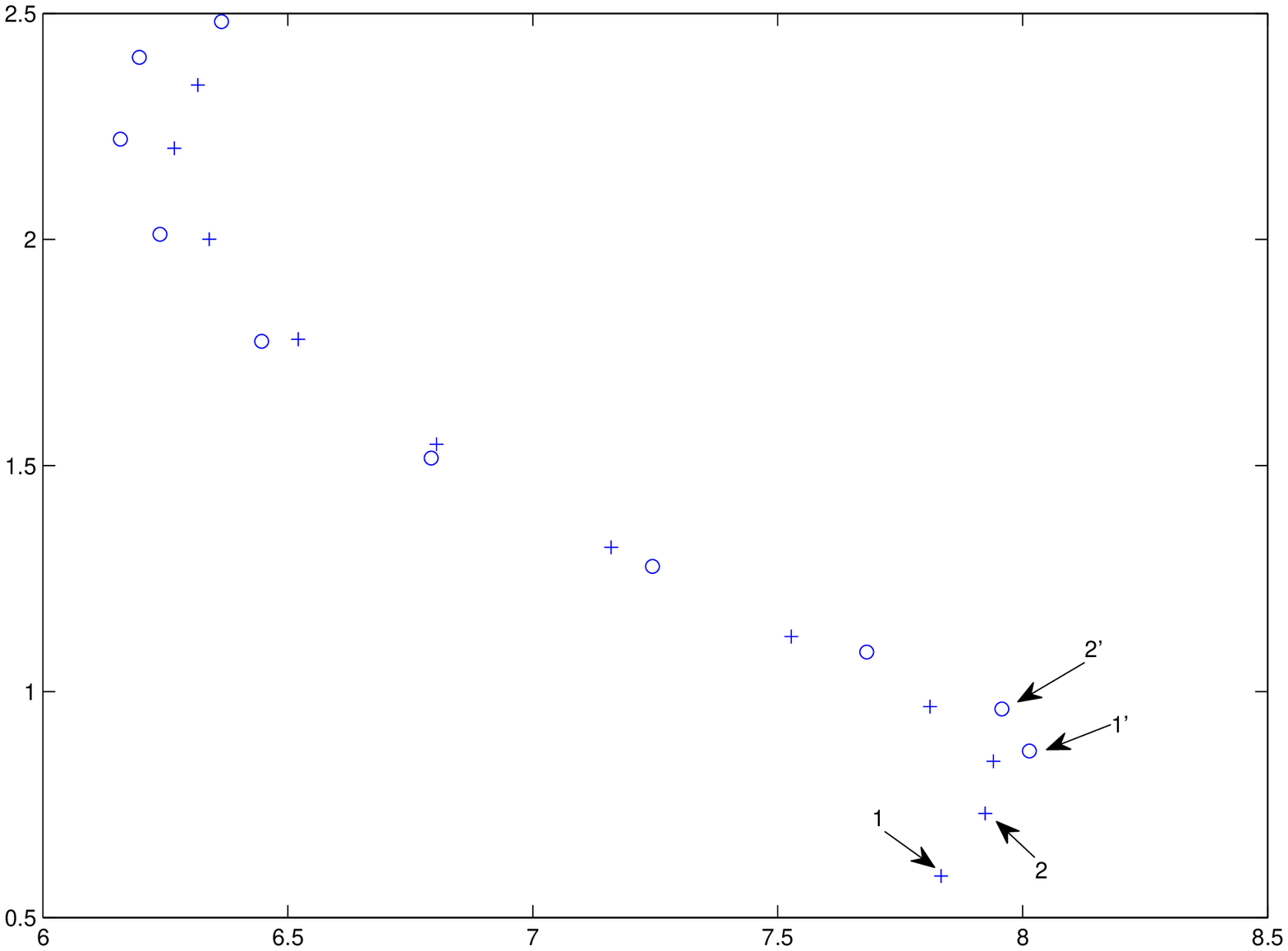,height=7.5cm,width=7.5cm}}}
\caption{ Figure illustrating the process of a robot, following
the respective centroid (a) with SDS and (b) with CDS. In both cases the robot location in each time step is shown
by `+', and corresponding centroid locations are shown by `o'.} \label{follow_centroid_cds}
 \centerline{
\subfigure[]{\psfig{figure=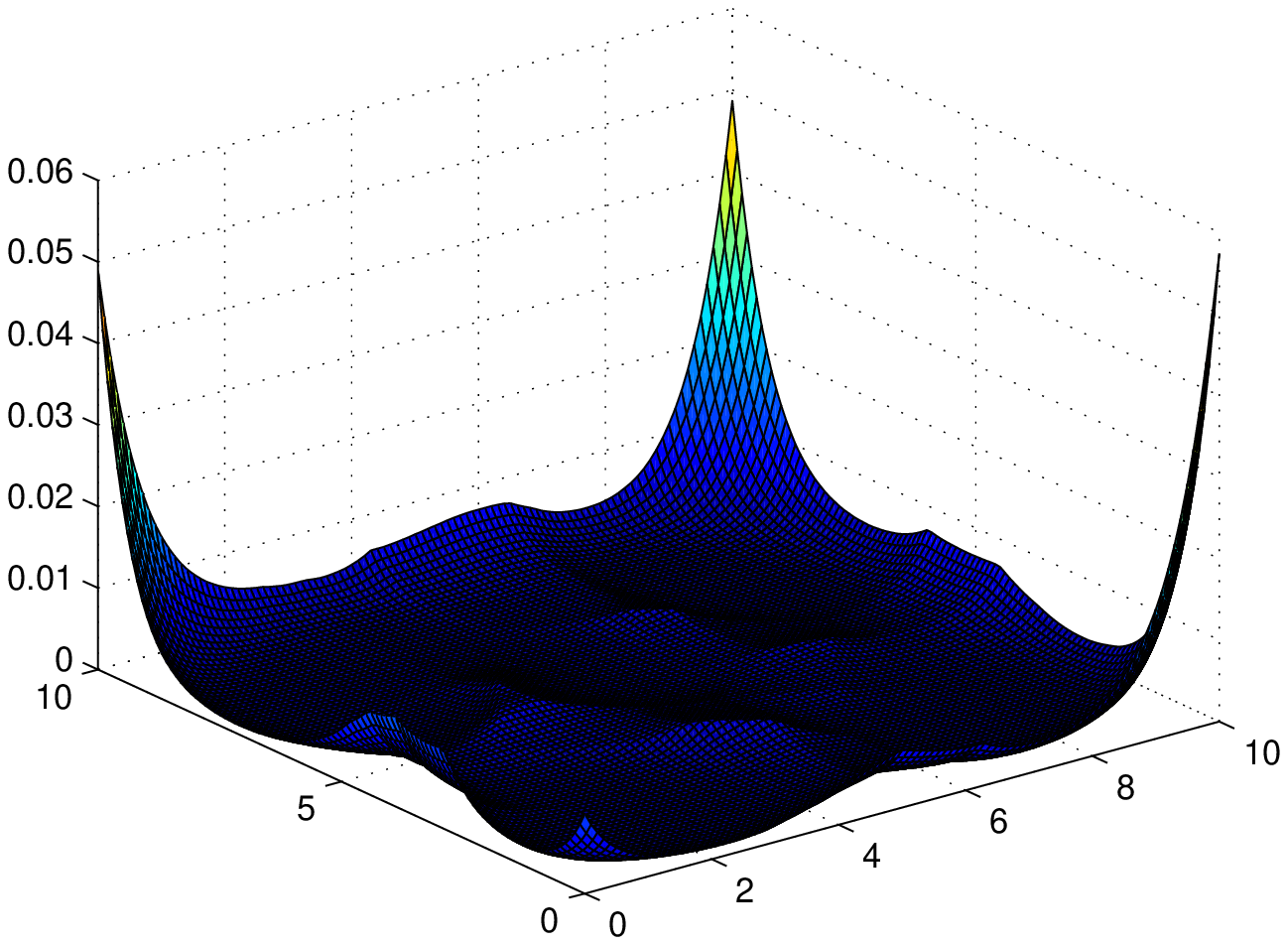,height=7.5cm,width=7.5cm}}
\subfigure[]{\psfig{figure=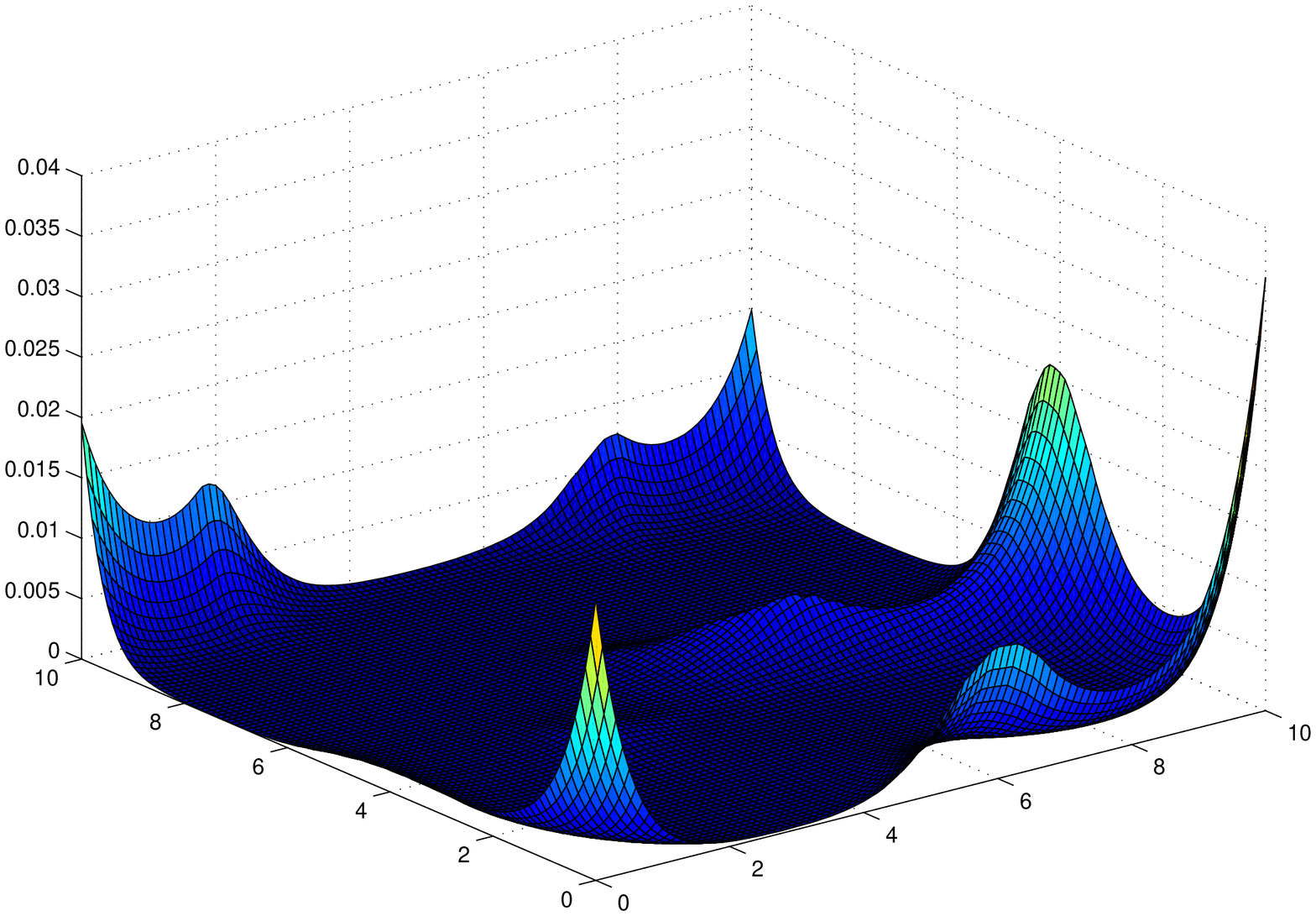,height=7.5cm,width=7.5cm}}
} \caption{The final uncertainty density distribution for (a) SDS
and (b) CDS strategies, with $N=5$ and without sensor range
limits} \label{dens_chap4}
\end{figure}

Figure \ref{follow_centroid_cds}(a) illustrates `robot 2'  moving toward centroid corresponding to its Voronoi cell with SDS strategy.
Robot's positions are marked with `+' while `o' marks the
centroids at successive time instances. Positions of robot in
first two time steps are marked as 1 and 2, while those of
centroids marked with $1'$ and $2'$. It can be observed that the robot
is tracking the centroid, which is changing as the Voronoi
cell is changing. Deployment stops and search is performed when
the robot is sufficiently close to the corresponding centroid.
One of the search instances is also marked, where, after search,
in order to track the next centroid, the robot takes an abrupt
turn. This leads to a non-smooth trajectory. Figure \ref{follow_centroid_cds} (b) illustrates one of the robots
moving toward centroid corresponding to its Voronoi cell with CDS strategy. It can be observed that the robot is tracking
the centroid, which is changing as the Voronoi cell and
the uncertainty density are changing.

\section{Implementation Issues} Here we discuss some of the
theoretical and implementation issues involved in  {\em combined
deploy and search} strategy.

\subsection{Spatial distributedness}
Here we discuss the implication of spatial distributedness of the proposed search strategies from a practical
point of view. We have seen that both the search strategies are spatially distributed in the Delaunay graph. These
results imply that all the robots need to do computations based on only local information, that is, by the knowledge
about position of neighboring robots. Also, the robots should have access to the updated uncertainty map within
their Voronoi cells. This can be achieved in several ways. One such way is that all the robots communicate with
a central information provider. But it is not necessary to have this global information. The $i$-th robot can communicate with its Voronoi neighbors (($\mathcal{N}_{\mathcal{G}}(i)$) and obtain the updated uncertainty information in a region $\cup_{\mathcal{N}_{\mathcal{G}}(i)} V_i$. As the Voronoi partition $\{V_i\}$ depends at least continuously on $\mathcal{P}$, the robot configuration \cite{bullo1}, in an evolving Delaunay graph, the communication within the neighbors is sufficient for each robot to obtain the uncertainty within its new Voronoi cell. The issues related
to communication of uncertainty information are not addressed in the paper except to assume that uncertainty
information is available to the robots. It is also possible that the robots can estimate the uncertainty map as done
in \cite{schwager}.

In practical conditions, the robots can communicate with other robots only when they are within the limits of
the sensor range. The Delaunay graph does not allow sensor range limitations to be incorporated. We need to use
 \emph{r-Delaunay} graph $\mathcal{G}_{LD}$  to incorporate the sensor range limitations. The scenario changes with incorporation of sensor range limitations into the search strategies. The updating of uncertainty density will also be within the sensor range limits (in fact, it is within the region $V_i \cap \bar{B}(p_i,R)$). The centroid that is computed will also be within the new restricted area. For an optimal deployment problem, from the perspective of sensor coverage, it has been observed that the corresponding control law is still spatially distributed (in r-limited Delaunay graph) and globally asymptotically stable.

\subsection{Synchronization} Synchronization plays an important role in multi-agent systems. Here we discuss this issue for both SDS and
CDS strategies. In the case of SDS, theoretically all the robots reach the respective centroid at infinite time. But in a
practical implementation, the robots are required to be sufficiently close, where the closeness is suitably defined, to
the respective centroids before starting the search operation. It is possible that at any point in time, different robots
are at different distances from the corresponding centroid. The robots need to come to a consensus as to when to
end the deployment and perform search operation. We have implemented the strategy in a single centralized program using MATLAB. In a practical situation, synchronization can be attained by robots communicating a flag bit indicating if a robot has reached its centroid or not. When all the robots  have reached the respective centroid within a tolerance distance, the search can be performed. We also assume that sensing and communication are instantaneous. In our simulation experiments we assume such a communication exists. Since the objective of this work is to evaluate the effectiveness of the search strategies, we make assumptions that simplify implementation without affecting the search effectiveness.

 CDS operates in a synchronous manner by design. If all the robots start at the same instant of time and have synchronized clocks, the search task is performed by every robot after the same interval of time. Given an accurate
global clock, synchronization is not a major issue in case of CDS.

Further in \cite{bullo1}, authors provide an asynchronous implementation for coverage control which can be suitably modified for CDS to operate asynchronously.

\section{Results and Discussion} In this section we present simulation experiment results to compare the performance of SDS and CDS strategies with {\em greedy search} and {\em random search} tailored to suit the problem setting addressed in this work. {\em Random} and {\em greedy} are generic strategies which can be tailored to suit any problem setting. We restrict the robots to move at a constant speed and also discretize the heading direction of robots at 1 degree resolution. Below, we discuss simple strategies used for the purpose of comparing the performance.

Typically search problems do not consider dynamics of search robots as the focus is more on the effectiveness
of search, that is, being able to identify region of high uncertainty and distribute search effort to reduce uncertainty.
Moreover, it is usually assumed that the search region is large compared to the physical size of the robot or the
area needed for the robot to maneuver. We assume a first order dynamics for the robots for the purpose of simulations.

\subsection*{Greedy search}

We discuss two simple {\em greedy} strategies. In the CDS, the $i$-th robot moves toward the centroid of
$V_{i}\cap\bar{B}(p_i,R))$ with $\tilde\phi$ as density. The Voronoi
partition takes care of the coordination among the neighboring
robots. In greedy search we simply let the $i$-th robot move toward
corresponding centroid of $\bar{B}(p_i,R)$. Thus there is no
cooperation between robots.

According to the update rule (\ref{phiupdate}), only the robot which
is most effective will perform search task at any given point $q \in
Q$. This leads to the idea of each robot performing search task
within its own Voronoi cell. We call such a greedy strategy as
{\em Voronoi greedy
 search} (VGS) strategy. In this strategy, only the control law is
greedy whereas the search is performed in a cooperative manner. In
case of cooperative search, each robot discards information about
the area which is better accessible to its neighboring robots.

In a true greedy search strategy we expect no cooperation among
robots at any point in time. A {\em true greedy strategy} (TGS) can
be achieved if each robot performs search task independently at
every point within $\bar{B}(p_i,R)$, leading to duplication of the
search task. The update law takes the form,
\begin{equation}
\label{phi_update_multiple} \phi_{n+1}(q) = \phi_n(q) \prod_{\{i |
p_i \in \bar{B}(q,R)\}} \beta(\parallel p_i - q
\parallel)
\end{equation}
It should be noted that the idea of collecting information by all
sensors in contrast to using the information from the most effective
sensor is problem dependent. If the sensor is a camera, then many
cameras taking picture of the same area may not add to the
information.

\subsection*{Random search}

Random search (RS) is probably the simplest search strategy. Here we
assume the robots move with a constant speed and the direction of
the robots is generated randomly. The uncertainty density updating
can be  given either by (\ref{phiupdate}) or by
(\ref{phi_update_multiple}). We have used the update law given by
(\ref{phiupdate}) in this case.

In TGS, VGS, and RS strategies, the search task is performed in every time step as in CDS.

\subsection*{Simulation results for performance comparison}The comparison
is based on the number of time steps required to reduce the average
uncertainty below a specified value.

The parameters which were varied during the simulations are the
number of robots ($N$), the sensor range limit ($R$), and the speed
of robots ($U$). It was desired that the average uncertainty density
should be reduced to a value below $0.002$.
\begin{table}
\begin{scriptsize}
\begin{center}
\begin{tabular}{|l|r|r|r|r|r|}
\hline
& CDS & TGS & VGS & SDS & RS \\
&&&& search steps  (time) &\\
\hline \hline
5.1.10 & 427 & 615 & 1115 & 193 (1207) & 2000+ \\
\hline
5.1.25 & 338 & 541 & 1416 & 192 (520) & 2000+ \\
\hline
5.1.50 & 1292 & 1150 & 1856 & 193 (261) & 1931 \\
\hline
5.2.10 & 264 & 285 & 250 & 53 (850) & 2000+ \\
\hline
5.2.25 & 99 & 217 & 251 & 53 (359) & 1188 \\
\hline
5.2.50 & 81 & 179 & 267 & 53 (164) & 814 \\
\hline
5.2.100 & 112 & 96 & 209 & 52 (82)& 246\\
\hline
5.3.10 & 100 & 202 & 538 & 25 (532)  & 2000+ \\
\hline
5.3.25 & 76 & 97 & 193 & 26 (239)& 1645 \\
\hline
5.3.50 & 42 & 66 & 177 & 25 (115)  & 427 \\
\hline
5.3.100 & 35 & 49 & 127& 25 (52) & 252\\
\hline
5.3.150 & 42 & 39 & 144 & 25 (34.7) & 46 \\
\hline
5.3.200 & 42 & 43 & 57 & 25 (26) & 54 \\
\hline
5.3.300 & 132 & 63 & 83 & 25 (17.34) & 56\\
\hline
5.4.10 & 61 & 61 & 73 & 17 (404)  & 2000+ \\
\hline
5.4.25 & 29 & 74 & 87 & 17 (141.6) & 509 \\
\hline
5.4.50 & 33 & 39 & 99 & 17 (77.25) & 230 \\
\hline
5.4.100 & 20 & 24 & 80& 17 (41.5) & 74 \\
\hline
5.4.150 & 26 & 43 & 216& 17 (27.6)& 52 \\
\hline
5.4.200 & 28 & 30 & 72 & 17 (20.74) &54 \\
\hline
5.4.400 & 34 & 36 & 47 & 17 (10.37) & 31 \\
\hline
20.1.10 & 124 & 264 & 429 & 48 (397.6)  & 2000+ \\
\hline
20.1.25 & 79 & 200 & 846 & 48 (163.5) & 2000+ \\
\hline
20.1.50 & 96 & 92 & 909 & 48 (81.7) & 2000+ \\
\hline
20.2.10 & 70 & 71 & 109 & 15 (217.6)  & 2000+ \\
\hline
20.2.25 & 31 & 43 & 65 & 15 (86.0) & 375 \\
\hline
20.2.50 & 23 & 31 & 123 & 15 (45.0) & 305 \\
\hline
20.2.100 & 25 & 23 & 135 & 15 (22.5) & 73 \\
\hline
20.2.150 & 33 & 22 & 56 & 15 (15.0) & 50 \\
\hline
20.3.10 & 61 & 58 & 64 & 9 (131.1)  & 2000+ \\
\hline
20.3.25 & 28 & 28 & 56 & 9 (52.6) & 398 \\
\hline
20.3.50 & 17 & 17 & 70 & 9 (25.8) & 254 \\
\hline
20.3.100 & 14 & 12 & 76 & 9 (12.8) & 48 \\
\hline
20.3.150 & 16 & 11 & 30 & 9 (8.5) & 29 \\
\hline
20.3.200 & 20 & 11 & 22 & 9 (6.4) & 26 \\
\hline
20.3.300 & 38 & 13 & 46 & 9 (4.2) & 20\\
\hline
20.4.10 & 54 & 50 & 58 & 8 (120)  & 1599 \\
\hline
20.4.25 & 25 & 24 & 31 & 8 (38.7) & 615 \\
\hline
20.4.50 & 15 & 15 & 94 & 8 (22.2) & 160 \\
\hline
20.4.100 & 11 & 10 & 77 & 8 (11.2) & 45 \\
\hline
20.4.150 & 11 & 10 & 98 & 8 (6.3) & 30 \\
\hline
20.4.200 & 13 & 12 & 47 & 8 (4.8) & 18 \\
\hline
20.4.400 & 95 & 9 & 30 & 8 (2.4) & 11\\
\hline \hline
\end{tabular}\label{comparison_table}
\end{center}
\end{scriptsize}
\caption{\footnotesize Comparison of performance of various
multi-robot search strategies (CDS: Combined deploy and search; SDS:
Simultaneous deploy and search; VGS: Voronoi greedy search; TGS:
True greedy search; RS: Random search. The first column gives the
parameter of simulation experiments as (N).(R).(100U). The numbers
indicate the number of search tasks performed before reaching the
termination condition. The maximum number of steps was restricted to
2000. The figures within parentheses for the case {\em sequential
deploy and search} are the actual times required for achieving the
task.}
\end{table}

Table I compares the performance of the search
strategies discussed in this paper in terms of the number of
searches performed before the termination condition is reached. The
first column gives the parameters $N$, $R$, and $100\times U$ as
$N.R.100U$. Thus, $20.4.10$ means $N=20$, $R=4$, and $U=0.1$. From
Table I, it can be observed that in most cases,
the {\em combined deploy and search} strategy performs better
followed by {\em true greedy} strategy, {\em Voronoi greedy}
strategy, and random search in order of degrading performance as
indicated by the number of time steps. In terms of the number of
searches being performed, {\em sequential deploy and search}
strategy performs better than {\em combined deploy and search}
strategy.

\begin{figure*}
\centerline{
\subfigure[]{\psfig{figure=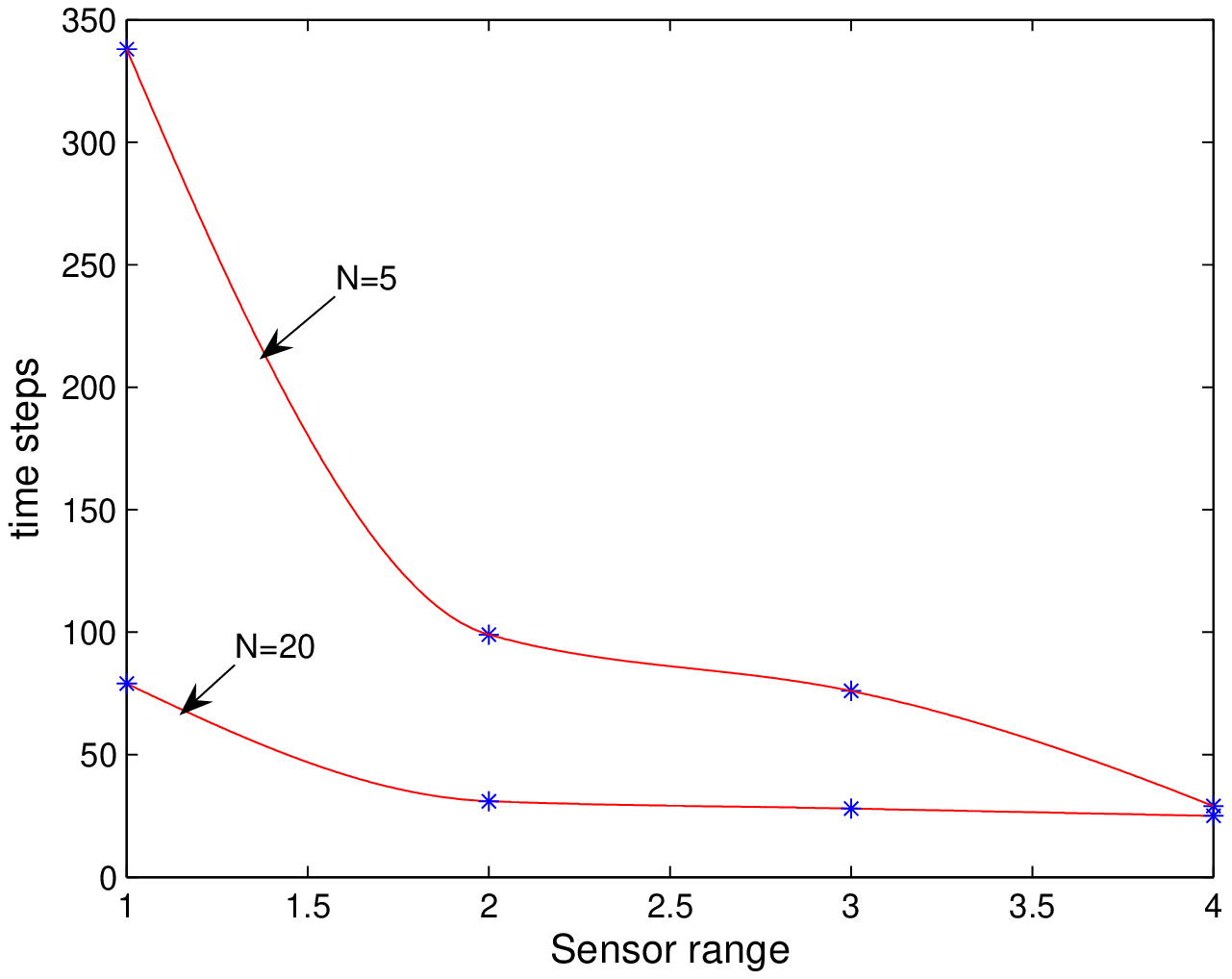,height=4.8cm,width=4.8cm}}
\hspace{-0.2in}
\subfigure[]{\psfig{figure=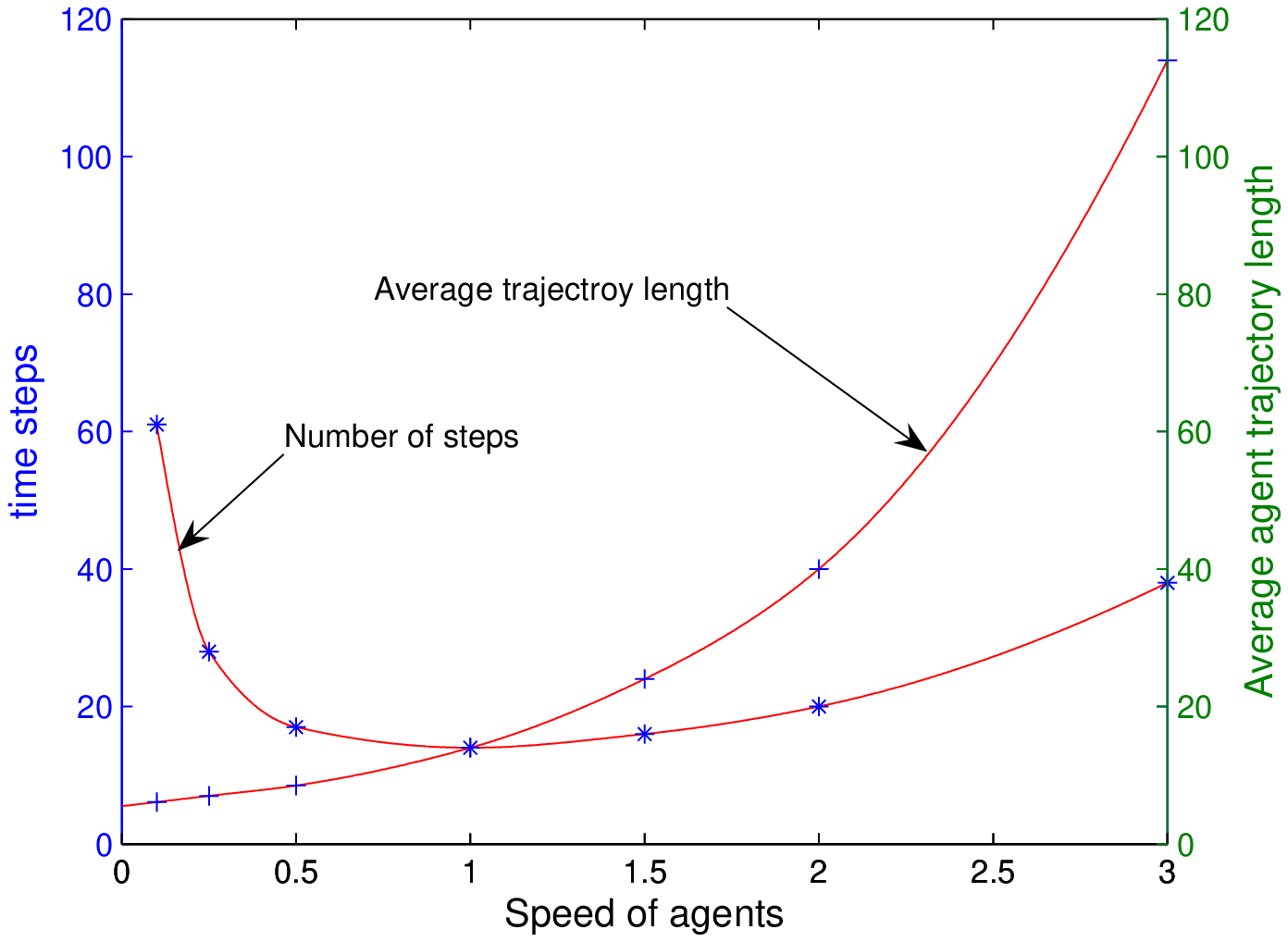,height=4.8cm,width=4.8cm}}
\hspace{-0.1in}
\subfigure[]{\psfig{figure=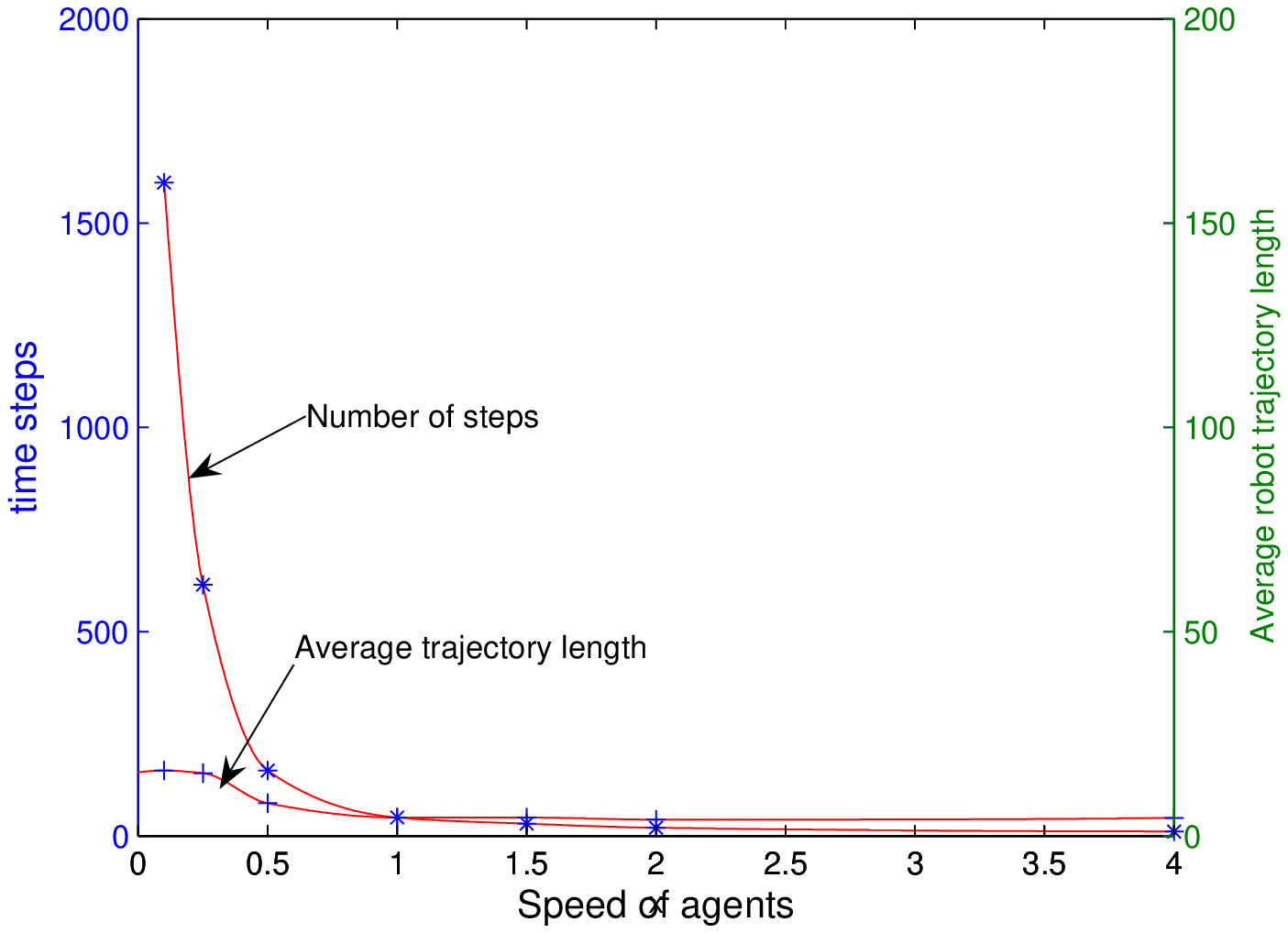,height=4.8cm,width=4.8cm}}}
\caption{\small (a) Number of time steps vs. sensor range $R$ for
{\em combined deploy and search} with a speed of $0.25$ for
different number of robots (b) Number of time steps versus robots'
speed $U$ for {\em combined deploy and search} for 20 robots with
sensor range $R =3$ (c) Number of time steps vs. robots' speed $U$
for {\em random search} for 20 robots with sensor range $R =4$.}
\label{cs_RU_rand_U}
\end{figure*}

The numbers in Table I indicate that with an
increase in $N$ and $R$ the performance of all the strategies
improves. Figure \ref{cs_RU_rand_U}(a) illustrates this. In this
figure, the points have been interpolated using a shape preserving
interpolation scheme.

Performance of CDS strategy deteriorates at higher speeds and
lower limit on sensor range as illustrated in Figure
\ref{cs_RU_rand_U}(b) for {\em combined deploy and search}. The robots do not move if the
centroid of respective $V_i\cap\bar{B}(p_i,R)$ is closer than
$U/2$, the tolerence set in the program. With large U and small $R$, quite often this condition
restricts the motion of robots. The condition adversely affects
the performance of CDS as the region $V_i\cap\bar{B}(p_i,R)$ is
always smaller than the region $\bar{B}(p_i,R)$, which is used by
greedy strategies for computing the centroid. That is, in case of
CDS, it is more likely that the centroids are closer than $U/2$ to
robots. In case of greedy strategies, the robots need to move
toward the centroids of corresponding $\bar{B}(p_i,R)$. This
region being relatively larger, it is less likely that the robots
will be closer to centroids. It has been observed during
simulations, particularly with CDS, that many of the robots do not
move during the entire search operation. When number of robots is
lower, difference in areas of above two regions is more likely
to be higher.

\begin{figure}
\centerline{
\subfigure[]{\psfig{figure=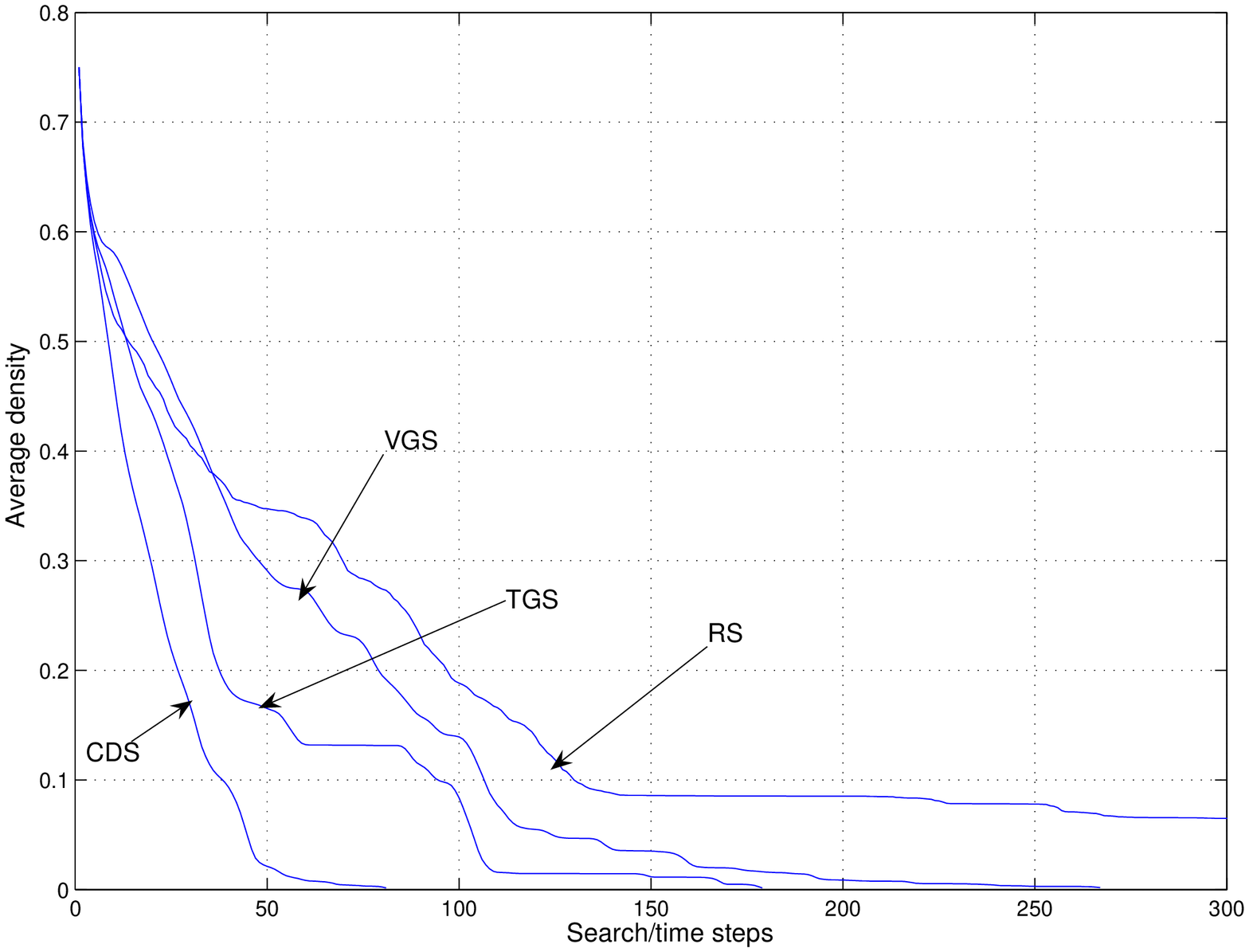,height=7.5cm,width=7.5cm}\label{avg_a}}
\subfigure[]{\psfig{figure=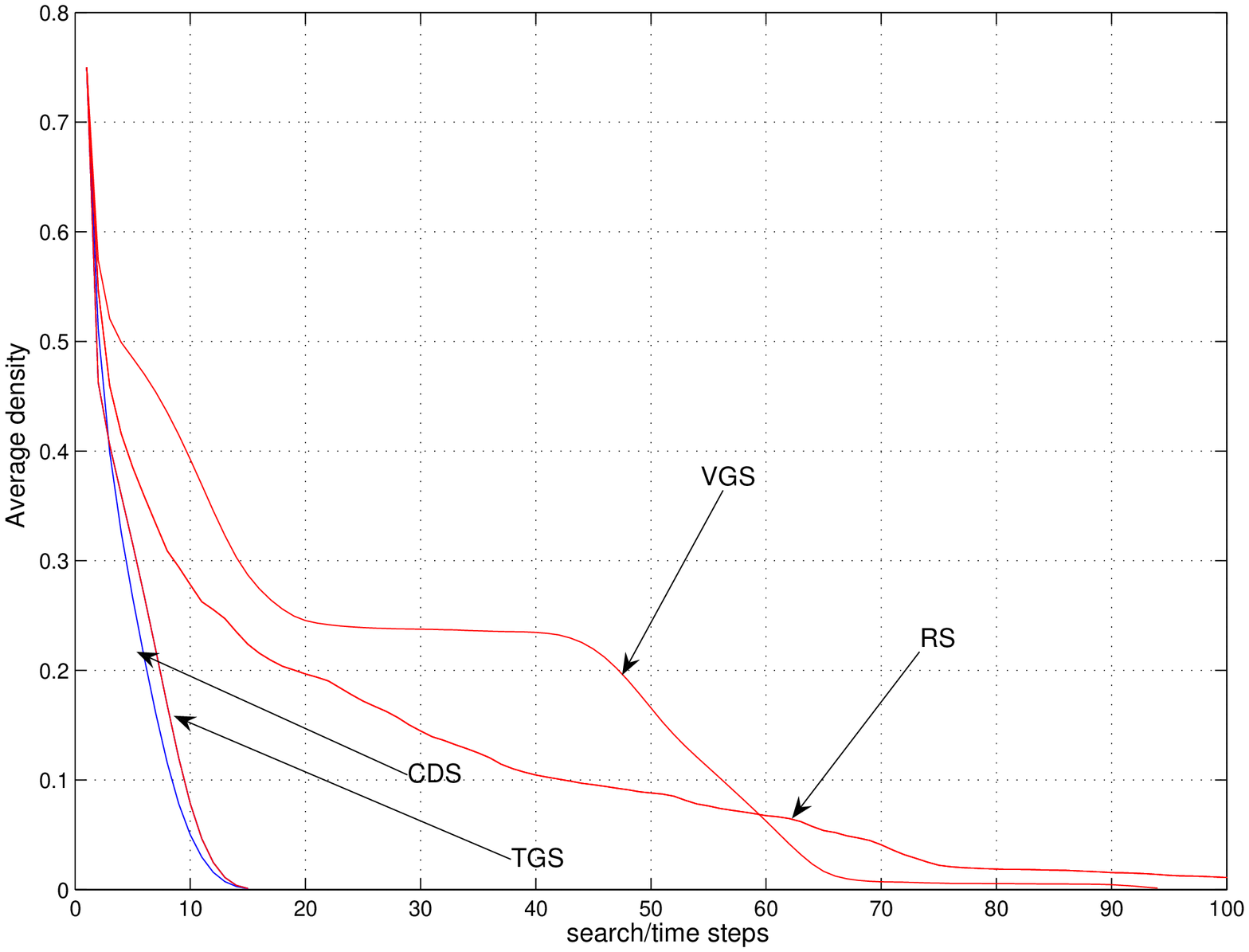,height=7.5cm,width=7.5cm}\label{avg_b}}
} \caption{Average uncertainty density distribution with $R=2$,
 $U=0.5$  and (a) with 5 robots and (b) with 20 robots.}
\label{avg_density}
\end{figure}

Figures \ref{avg_density} (a) and (b) show the time history of the
average uncertainty density distribution for {\em combined deploy
and search}, {\em true greedy}, {\em  Voronoi greedy}, and {\em
random} search strategies for two different sets of parameters, one
with $N=5$, $R=2$ units, and $U=0.5$ units and another with $N=20$,
$R=4$ units, and $U=0.5$ units. Both the figures illustrate that the
{\em combined deploy and search} takes minimum time steps while {\em
true greedy}, {\em  Voronoi greedy}, and {\em random} search
strategies follow in order of degrading performance in terms of
number of time steps. It should be noted though all these strategies perform search task in every time step, CDS performs better because of better coverage of the search space due to cooperation among the robots through the Voronoi partition. The CDS reduces uncertainty, in most cases, more than the other strategies (TGS, VGS, and RS) for same number of searches.

\begin{figure}
\centerline{
\subfigure[]{\psfig{figure=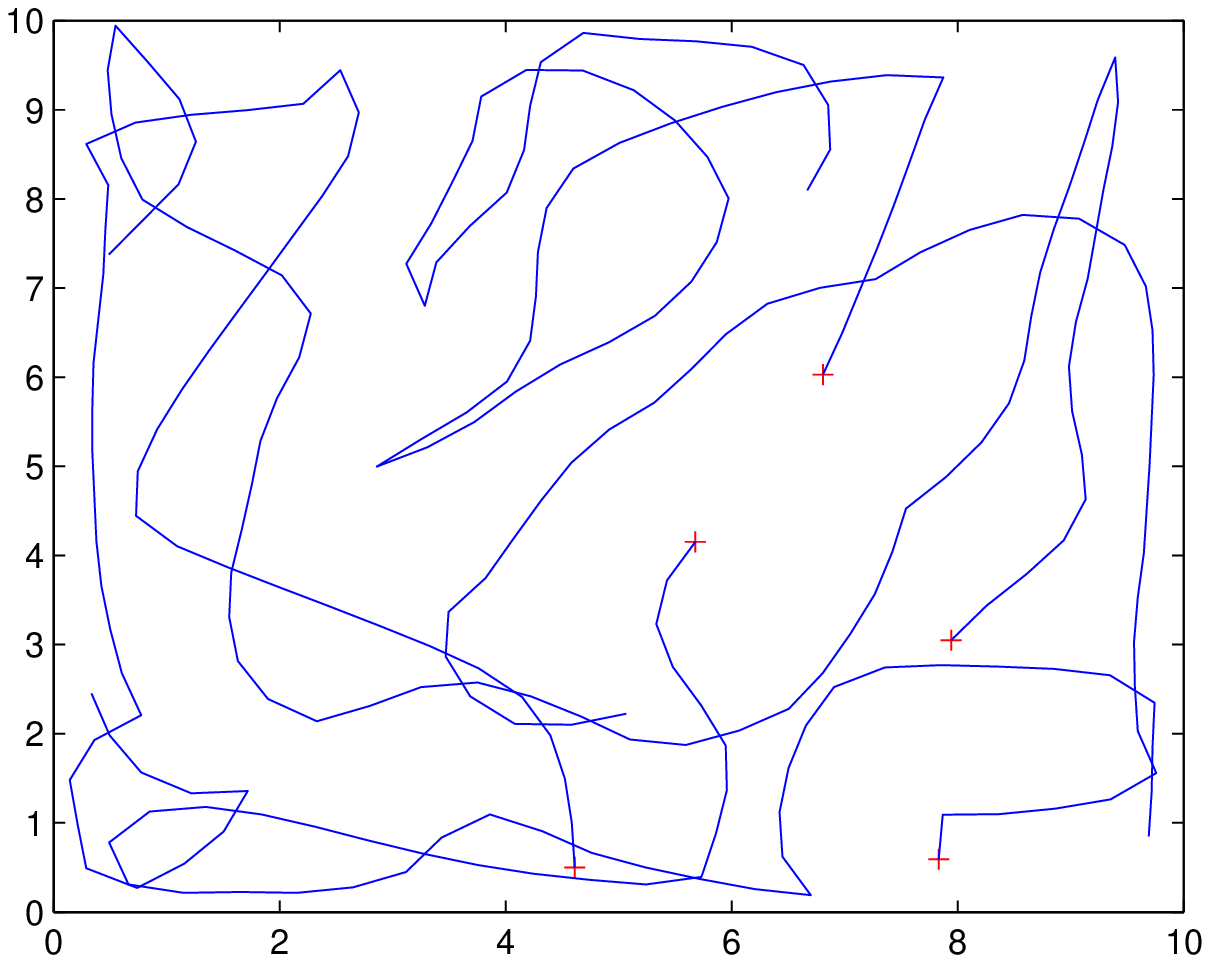,height=7.5cm,width=7.5cm}\label{cs_traj_5}}
\subfigure[]{\psfig{figure=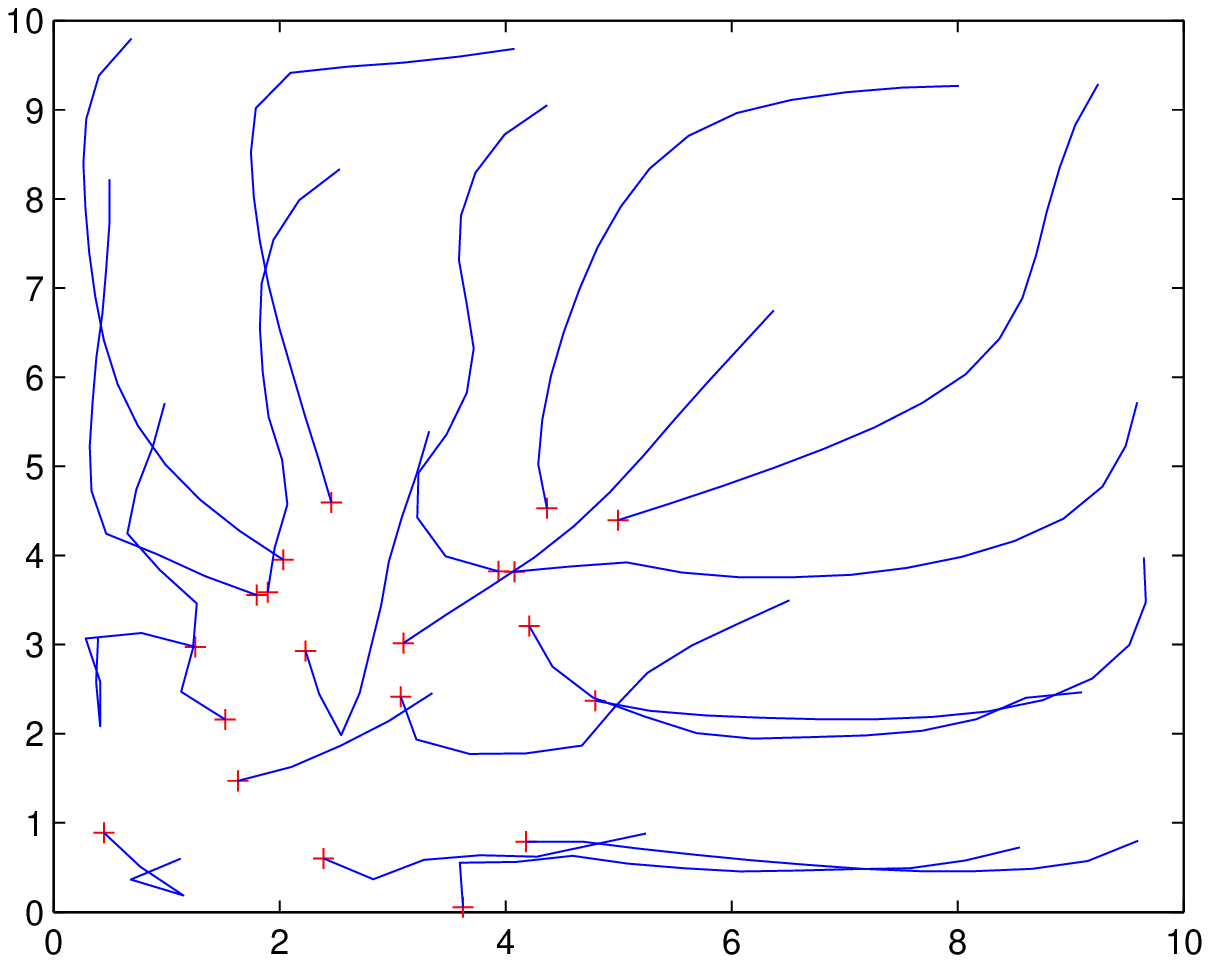,height=7.5cm,width=7.5cm}\label{cs_traj_20}}
} \caption{Trajectories of robots with {\em combined deploy and
search} strategy for  $U=0.5$  and (a) 5 robots with sensor range
$R=2$, (b) 20 robots with $R=4$. '+', indicate initial positions of
robots.} \label{cs_traj}
\end{figure}

\begin{figure}
\centerline{
\subfigure[]{\psfig{figure=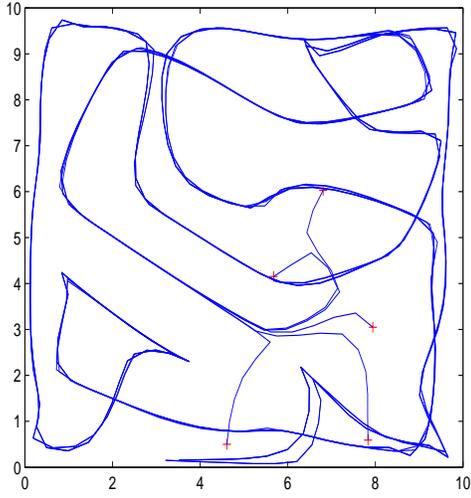,height=7.5cm,width=7.5cm}\label{greedy1_traj_a}}
\subfigure[]{\psfig{figure=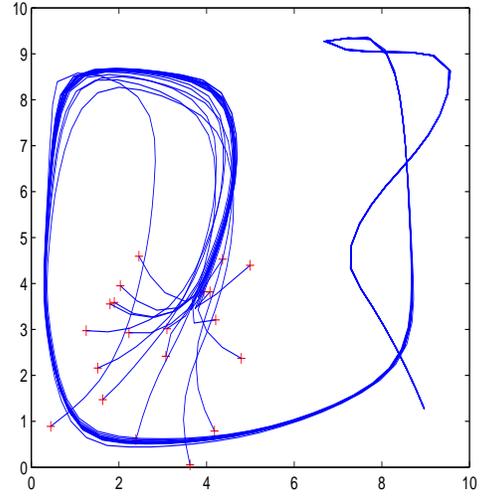,height=7.5cm,width=7.5cm}\label{greedy1_traj_b}
}} \caption{Trajectories of robots with {\em Voronoi greedy}
strategy for (a) $N=5$, $R=2$ and $U=0.5$ (b)$N=20$, $R=4$ and
$U=0.5$. '+', indicate initial positions of robots.}
\label{greedy1_traj}
\end{figure}

\begin{figure}
\centerline{
\subfigure[]{\psfig{figure=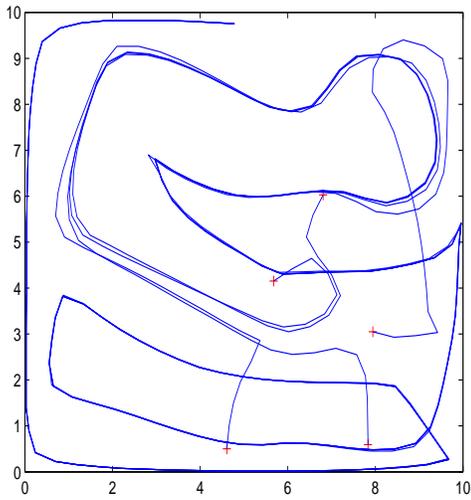,height=7.5cm,width=7.5cm}\label{greedyN_traj_a}}
\subfigure[]{\psfig{figure=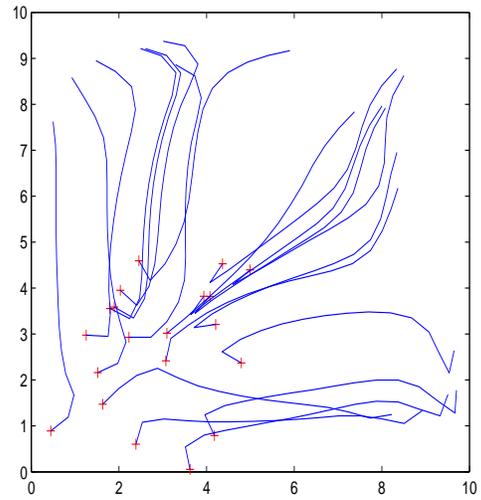,height=7.5cm,width=7.5cm}\label{greedyN_traj_b}}
} \caption{Trajectories of robots with {\em true greedy search}
strategy for (a) $N=5$, $R=2$ and $U=0.5$ (b)$N=20$, $R=4$ and
$U=0.5$. '+', indicate initial positions of robots.}
\label{greedyN5_2_50Traj}
\end{figure}

\begin{figure}
\centerline{
\subfigure[]{\psfig{figure=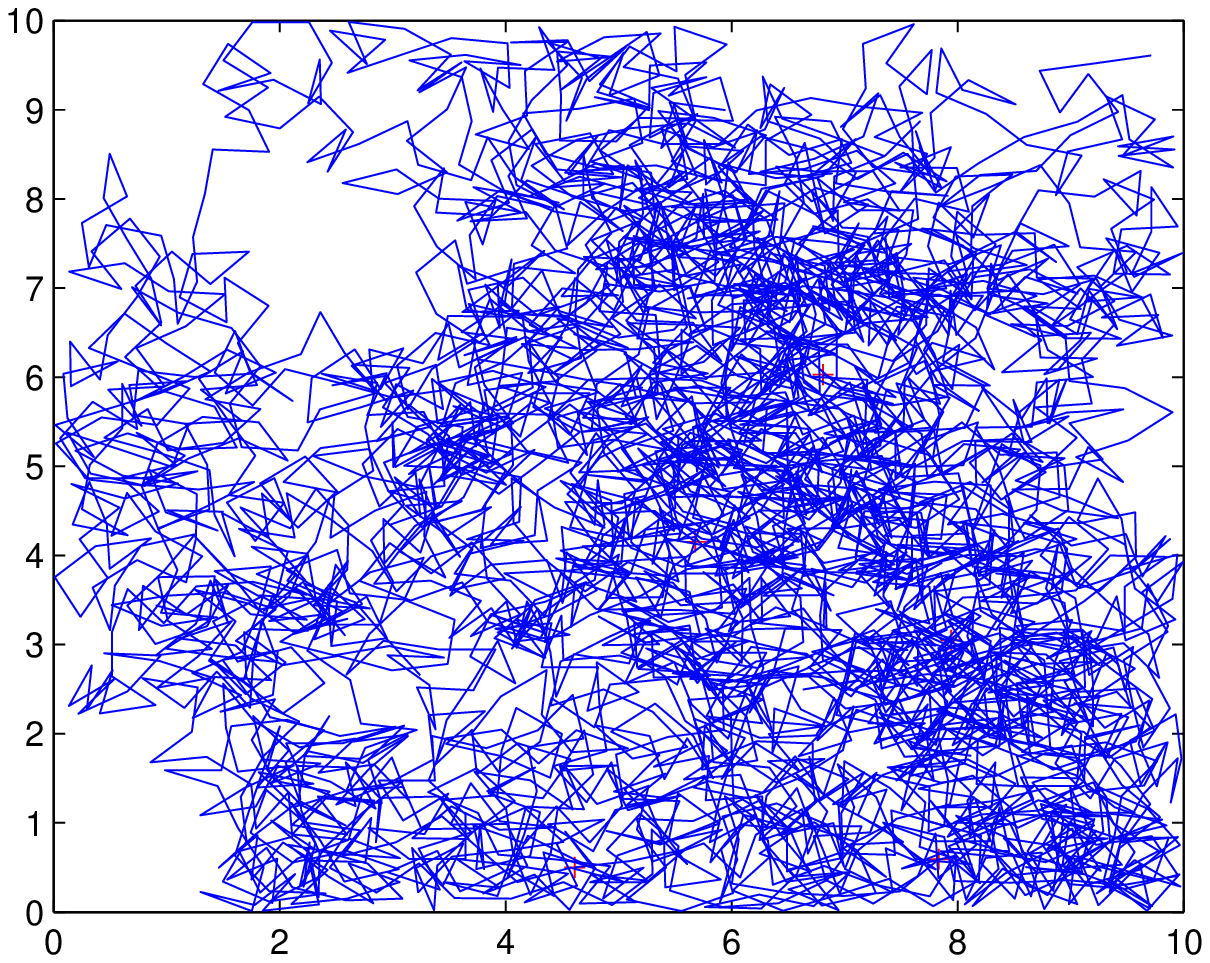,height=7.5cm,width=7.5cm}\label{random_traj_a}}
\subfigure[]{\psfig{figure=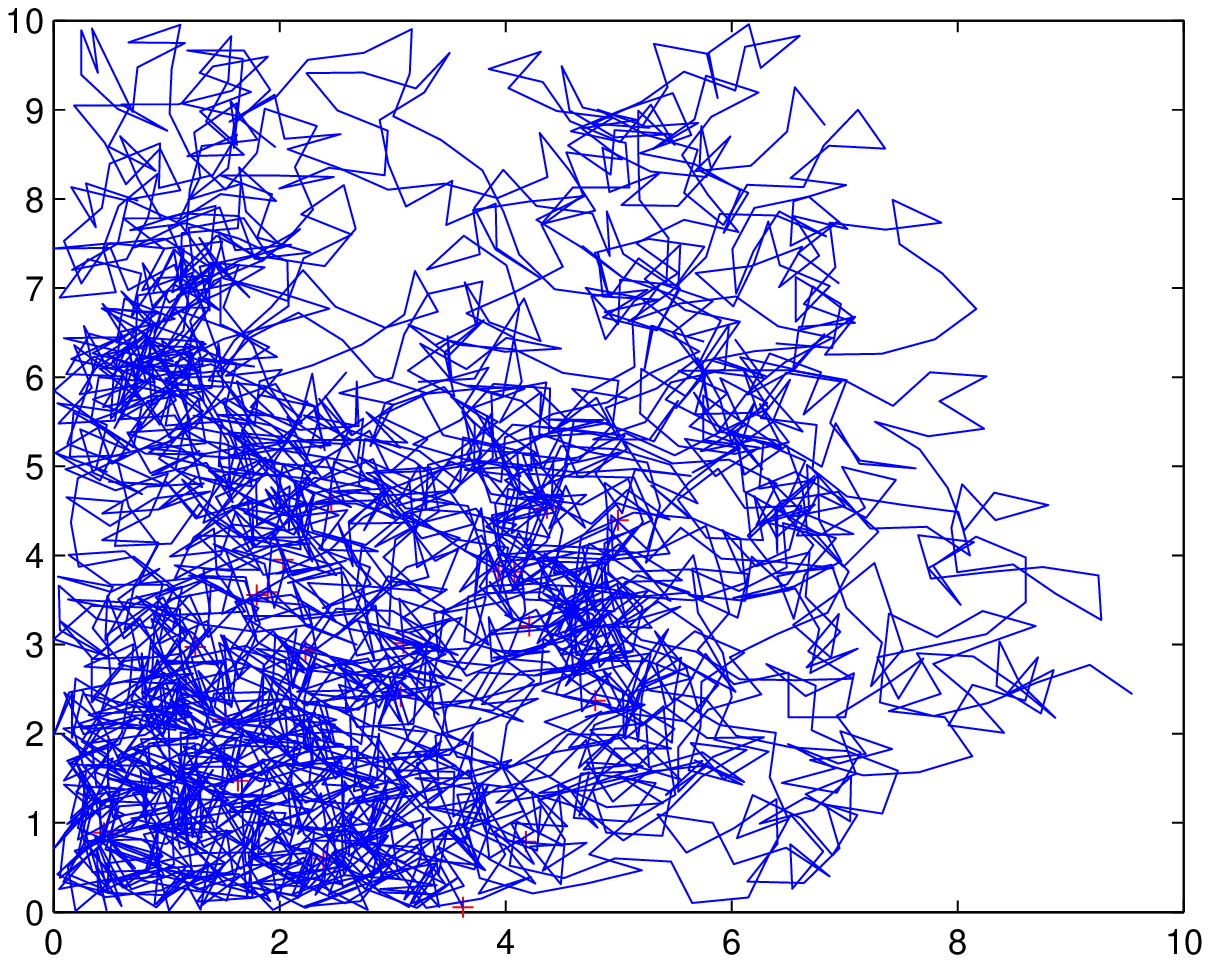,height=7.5cm,width=7.5cm}\label{random_traj_b}}
} \caption{Trajectories of robots with {\em random search} strategy
for (a) $N=5$, $R=2$ and $U=0.5$ (b)$N=20$, $R=4$ and $U=0.5$. The
initial positions of robots are as in Figure \ref{cs_traj}. The
number of points in each trajectory is 814.} \label{random_traj}
\end{figure}

Figures \ref{cs_traj} (a) and (b) show the trajectories of robots
for the above sets of parameters. Figures
\ref{greedy1_traj}-\ref{random_traj} show the robot trajectories
with {\em Voronoi greedy}, {\em true greedy}, and {\em random}
search strategies, respectively, for the same sets of parameters.

With the greedy strategies, when the parameters are such that the
convergence is slower, it is more likely that robots start moving
together in clusters. This happens due to lack of cooperation
between robots and each of them moving toward points with higher
uncertainty density. Once this happens, the robots move together,
and in case of {\em Voronoi greedy} search where only one of the
robot is active at a given point, the basic purpose of deploying
multiple robots for search is defeated, leading to a much slower
convergence. The strategy is equivalent to a single robot searching
the space. With {\em true greedy} search strategy, though all the
robots perform search within their sensor range, only in a very few
cases, the performance is comparable to {\em combined deploy and
search}. Cooperation among robots in {\em combined deploy and
search} and {\em sequential deploy and search} strategies ensure a
better coverage of the space and sharing of search load among
robots.

As expected, {\em random search} does not perform well when the
speed is less. But as the robot speed is increased, the {\em random
search} performs better. Unlike in other strategies, in case of the
{\em random search} strategy, the performance always increases with
speed as illustrated in Figure \ref{cs_RU_rand_U}(c). We compute the
average length of the trajectory  of a robot by multiplying the
time, which is equal to the number of steps with the speed. Figure
\ref{cs_RU_rand_U}(c) also shows the average trajectory length of
robots indicating that after a speed of 1 unit, average trajectory
length is almost a constant value. It is interesting to note that
the {\em random search} strategy performs better than all other
strategies at high robot speeds. This is due to better coverage of
the search space by robots. Motion of robots are restricted by the
specified control law in case of other strategies, which does not
happen in case of {\em random search}.

\begin{figure}
\centerline{
\subfigure[]{\psfig{figure=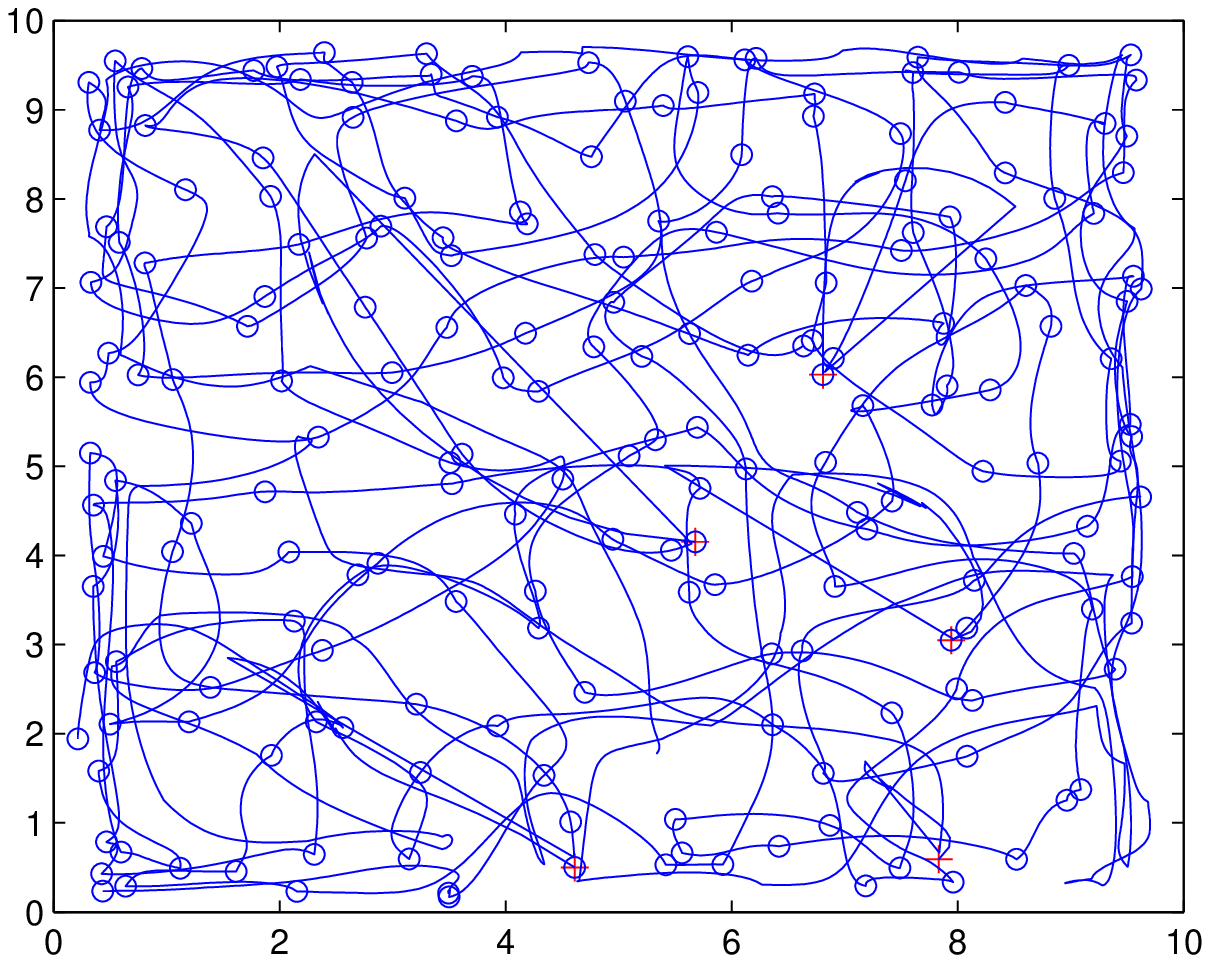,height=7.5cm,width=7.5cm}\label{dns_traj_a}}
\subfigure[]{\psfig{figure=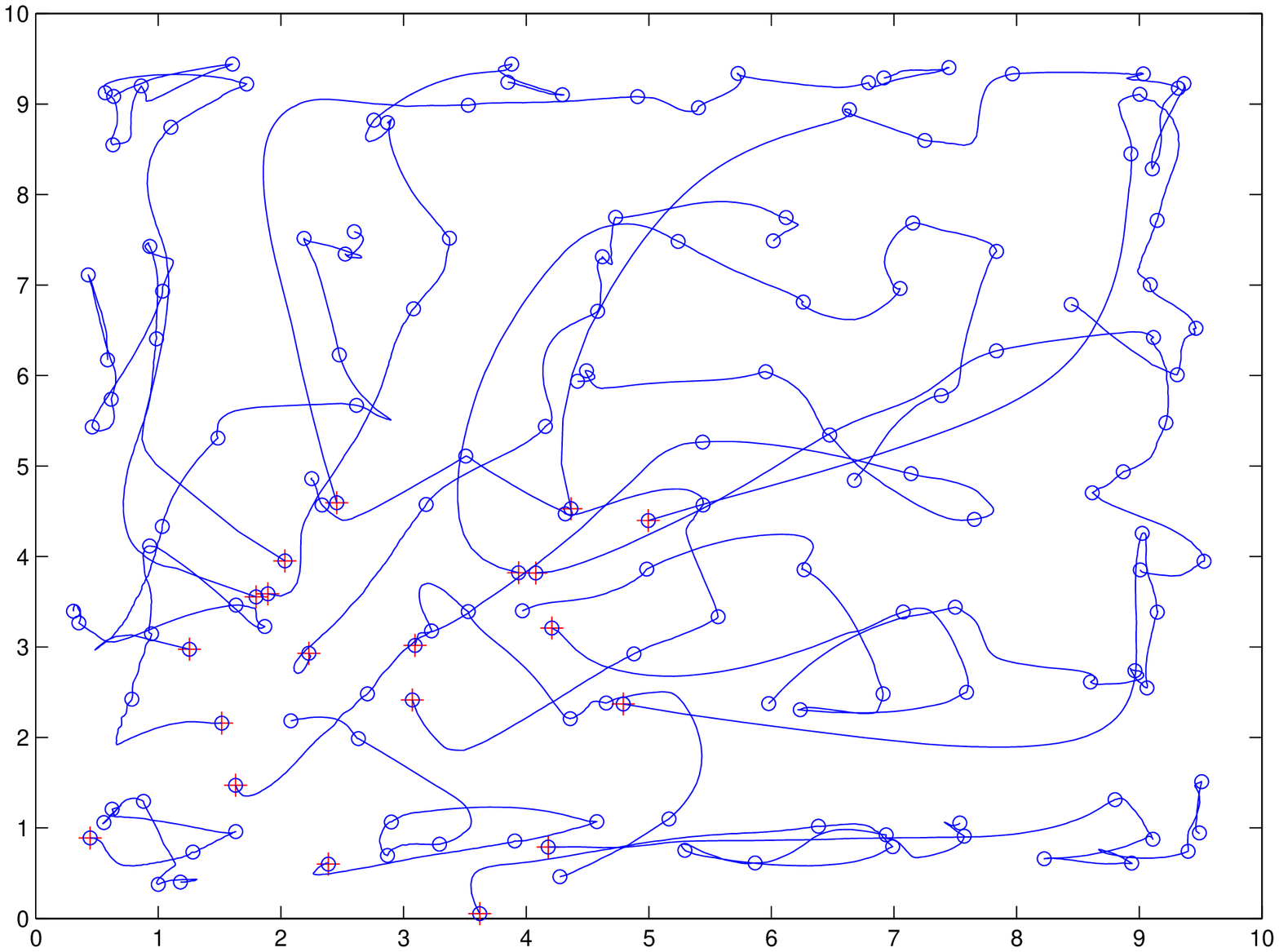,height=7.5cm,width=7.5cm}\label{dns_traj_b}}
} \caption{Trajectories of robots with {\em sequential deploy and
search} strategy for (a) $N=5$, $R=2$ and $U=0.5$ (b)$N=20$, $R=4$
and $U=0.5$. '+', indicate initial positions of robots and 'o' indicate end of deployment and points along the trajectories where search is performed.}
\label{dns_traj}
\end{figure}

In the case of {\em sequential deploy and search} strategy, as the
robots approach the optimal deployment, we allow partial stepping to
enable robots to move to as close as possible the respective
centroids. Figure \ref{dns_traj} shows robot trajectories for the
{\em sequential deploy and search} strategy.  With a
small sensor range, the time required for completing search is more
than that for the rest. But as the sensor range is increased, there
is a continuous improvement in the performance leading to the fastest
convergence amongst the compared strategies. Note that in the case
of the rest of the strategies, the time taken is the same as the
number of search steps. The performance of {\em sequential deploy
and search} strategy seems to monotonically increase with $N$, $R$
and $U$.

In our formulation of the objective function, we do not attach any
cost to the search operation. This is the basic motivation for
strategies such as {\em combined deploy and search}. In such a
situation we observe that the {\em combined deploy and search} is
best suited for reducing uncertainty. But in a practical situation,
some cost may be associated with each search operation, in which
case {\em sequential deploy and search} strategy is more efficient.
{\em True greedy} is effective in certain parameter ranges and when
in a situation where the duplication of information by more than one
sensor is beneficial. {\em  Voronoi greedy} is always below  both
{\em combined deploy and search} and {\em true greedy} in
performance in terms of total time.

It can be noted that, though the Voronoi partition based strategies
proposed in this paper are computationally complex, their
performance is better than the simple strategies such as {\em true
greedy} and {\em random search}. The Voronoi based strategies result
in collision free trajectories in a natural way, whereas non-Voronoi
based search strategies need additional collision avoidance schemes
and require additional computation.

\section{Conclusion and future work} The problem of multi-robot search in an unknown environment with a
known uncertainty probability distribution function is addressed in
this paper. A multi-robot search strategy namely  {\em combined deploy and search} (CDS) was presented as a modification of {\em sequential deploy and search} (SDS) strategy proposed in our previous work. We have shown that the centroidal Voronoi configuration
with respect to the density as perceived by the sensors are the instantaneous
critical points of the objective function maximizing the single step
search effectiveness, if the robots are located at corresponding centroids. It was shown, that the CDS strategy can reduce the average uncertainty to any arbitrarily small value in finite time. The optimal deployment strategy has been analyzed in presence of some constraints on robot speed and limit on sensor range for convergence of the robot trajectories with the corresponding control laws responsible for the motion of robots. The proposed CDS strategy was explained with help of a illustrative example, in comparison with SDS strategy.

Based on a set of simulation experiments, the performance of SDS and CDS strategies have been compared with standard strategies such as {\em greedy} and {\em random} search with a constant speed constraint on robots and sensor range limit. The parameters $N$, the
number of robots, $R$, the sensor range limit and $U$, the speed of
robots were varied to compare the performance of different
strategies. The simulation results indicated that the {\em
combined deploy and search} strategy is best suited when search
task does not involve any cost. The {\em true greedy} search
strategy performed well in a few parameter ranges and is acceptable
only in a situation where gathering of information by more than one
sensor in an area makes sense. When a cost is associated with the
search task, {\em sequential deploy and search} strategy was found to be more
suitable.

Though the Voronoi partition based strategies such as SDS and CDS
are computationally complex, their performance is better than the
simple non-Voronoi based strategies such as {\em true greedy} and {\em random search}.
The Voronoi based strategies result in collision free trajectories
in a natural way, whereas non-Voronoi based search strategies need
additional collision avoidance schemes and require additional
computation.

\pagebreak

\appendix
\begin{emph} Theorem A.1:  The gradient of the multi-center objective
function (\ref{obj1}) with respect to $p_i$ is given by
\begin{equation}
\label{grad} \frac{\partial \mathcal{H}_n}{\partial p_i} =
\int_{V_i}\phi(q)\frac{\partial}{\partial p_i}(1-\beta(r_i))dQ
\end{equation}
where $r_i = \parallel p_i - q \parallel$.
\end{emph}

\noindent {\it Proof.} Rewrite (\ref{obj1}) as
\begin{equation}
\mathcal{H} = \sum_{i\in\{1,2,\ldots,N\}} \mathcal{H}^i
\end{equation}
where $\mathcal{H}^i = \int_{V_i}(f(r))\phi(q)dQ$. Now,
\begin{equation}
\frac{\partial \mathcal{H}_n}{\partial p_i} = \sum_{j\in\{1,2,\ldots,N\}}
\frac{\partial \mathcal{H}_n^j}{\partial p_i}
\end{equation}
where, $f(\cdot) = 1-\beta(\cdot)$.

Applying the general form of the Leibniz theorem (\cite{kundu})
\begin{equation}
\label{grad_proof}
\begin{array}{lcl} \frac{\partial \mathcal{H}_n}{\partial p_i} &=&
\int_{V_i}\phi(q)\frac{\partial f}{\partial p_i}(\|q-p_i\|)dQ\\
&& + \sum_{j \in
N_i}\int_{A_{ij}}\mathbf{n}_{ij}(q).\mathbf{u}_{ij}(q)\phi(q)f(\|q-p_i\|)dQ\\
&&+ \sum_{j \in
N_i}\int_{A_{ji}}\mathbf{n}_{ji}(q).\mathbf{u}_{ji}(q)\phi(q)f(\|q-p_j\|)dQ
\end{array}
\end{equation}
where,
\begin{enumerate}
\item $N_i$ is the set of indices of agents which are neighbors of
the $i$-th  agent in the Delaunay graph $\mathcal{G}_{LD}$.

\item $A_{ij}$ is the part of the bounding surface (line segment in two dimensional case)
common to $V_i$ and $V_j$.

\item $\mathbf{n}_{ij}(q)$ is the unit outward normal to $A_{ij}$ at
$q \in A_{ij}$. Note that $\mathbf{n}_{ij}(q) =
-\mathbf{n}_{ji}(q)$, $\forall q \in A_{ij}$.

\item $\mathbf{u}_{ij}(q) = \frac{dA_{ij}}{dp_i}(q)$, the rate of
movement of the boundary at $q \in A_{ij}$ with respect to $p_i$.
Note that $\mathbf{u}_{ij}(q) = \mathbf{u}_{ji}(q)$.

\item Note also that $f(\|q-p_i\|) = f(\|q-p_j\|)$, $\forall q \in
A_{ij}$, as $\|q-p_i\|=\|q-p_j\|$, by definition of the Voronoi
partition.
\end{enumerate}

By (3)-(5) above, it is clear that for each $j \in N_i$,
\begin{displaymath}
\int_{A_{ij}}\mathbf{n}_{ij}(q).\mathbf{u}_{ij}(q)\phi(q)f(\|q-p_i\|)dQ
=
-\int_{A_{ji}}\mathbf{n}_{ji}(q).\mathbf{u}_{ji}(q)\phi(q)f(\|q-p_j\|)dQ
\end{displaymath}
Hence,
\begin{displaymath}
\frac{\partial \mathcal{H}_n}{\partial p_i} =
\int_{V_i}\phi(q)\frac{\partial f}{\partial p_i}(\|q-p_i\|)dQ
\end{displaymath}\hfill $\Box$

It can be observed that necessary smoothness conditions are valid as
\begin{enumerate}
\item[(i)] $\phi \text{, }f \in C^0$.

\item[(ii)] At the boundaries of the Voronoi cells, the
objective function does not have any jumps as $f(r_i) = f_j$ at the
boundary, by definition of the generalized Voronoi partitions.

\item[(iii)] The Voronoi partition itself varies smoothly  with $\mathcal{P}$.
\end{enumerate}
\end{document}